\RequirePackage{ifpdf}
\ifpdf % We are running pdfTeX in pdf mode
\documentclass[pdftex]{sigma}
\else
\documentclass{sigma}
\fi

\numberwithin{equation}{section}
\numberwithin{theorem}{section}
\numberwithin{proposition}{section}
\numberwithin{lemma}{section}
\numberwithin{corollary}{section}
\numberwithin{definition}{section}
\numberwithin{example}{section}
\numberwithin{remark}{section}
\numberwithin{note}{section}

\newtheorem{openq}{Open Question}[section]

\DeclareMathOperator{\Dom}{Dom}   %% domain of an operator
\DeclareMathOperator{\Tr}{Tr}     %% operator trace
\DeclareMathOperator{\tr}{tr}     %% matrix trace

\newcommand{\nm}[1]{\mbox{\ensuremath{\| #1 \|}}}
\newcommand{\fa}{\ensuremath{\forall}}

\newcommand{\inset}[2]{\ensuremath{\{ #1 \, | \, #2 \} }}
\newcommand{\inprod}[2]{\ensuremath{\langle #1 , #2 \rangle}}

\newcommand{\ceil}[1]{\ensuremath{\left\lceil #1 \right\rceil}}
\newcommand{\floor}[1]{\ensuremath{\left\lfloor #1 \right\rfloor}}

\newcommand{\A}{\ensuremath{\mathcal{A}}}

\newcommand{\ds}{\ensuremath{\langle D \rangle}}

\newcommand{\CC}{\ensuremath{\mathbb{C}}}
\newcommand{\RR}{\ensuremath{\mathbb{R}}}

\newcommand{\TT}{\ensuremath{\mathbb{T}}}
\newcommand{\SB}{\ensuremath{\mathbb{S}}}
\newcommand{\ZZ}{\ensuremath{\mathbb{Z}}}

\newcommand{\NN}{\ensuremath{\mathbb{N}}}

\newcounter{prop2count}

\begin{document}

\allowdisplaybreaks

\renewcommand{\thefootnote}{$\star$}

\renewcommand{\PaperNumber}{072}

\FirstPageHeading

\ShortArticleName{Measure Theory in Noncommutative Spaces}

\ArticleName{Measure Theory in Noncommutative Spaces\footnote{This paper is a
contribution to the Special Issue ``Noncommutative Spaces and Fields''. The
full collection is available at
\href{http://www.emis.de/journals/SIGMA/noncommutative.html}{http://www.emis.de/journals/SIGMA/noncommutative.html}}}

\Author{Steven LORD~$^\dag$ and Fedor SUKOCHEV~$^\ddag$}

\AuthorNameForHeading{S.~Lord and F.~Sukochev}

\Address{$^\dag$~School of Mathematical Sciences, University of Adelaide, Adelaide, 5005, Australia}
\EmailD{\href{mailto:steven.lord@adelaide.edu.au}{steven.lord@adelaide.edu.au}}

\Address{$^\ddag$~School of Mathematics and Statistics, University of New South Wales, Sydney, 2052, Australia}
\EmailD{\href{mailto:f.sukochev@unsw.edu.au}{f.sukochev@unsw.edu.au}}

\ArticleDates{Received March 25, 2010, in f\/inal form August 04, 2010;  Published online September 16, 2010}

\Abstract{The integral in noncommutative geometry (NCG) involves a non-standard trace called a Dixmier trace.  The geometric origins
of this integral are well known.  From a~measure-theoretic view, however, the formulation contains several dif\/f\/iculties.
We review results concerning the technical features of the integral
in NCG and some outstanding problems in this area.  The review is aimed
for the general user of NCG.}


\Classification{46L51; 47B10; 58B34}

%\vspace{-5mm}\tableofcontents

{\small \tableofcontents}

\renewcommand{\thefootnote}{\arabic{footnote}}
\setcounter{footnote}{0}

\section{Introduction} \label{sec:intro}

The basic object of noncommutative geometry (NCG) is a spectral triple
$(\A,H,D)$ \cite{CC}, originally called a K-cycle~\cite[p.~546]{CN}.
Here $\A$ is a countably generated unital
non-degenerate
$*$-algebra of bounded linear operators on a separable complex Hilbert space $H$,
and $D : \Dom D \to H$ is a selfadjoint linear operator with the properties $\nm{[D,a]} < \infty$
for all $a \in \A$ and $\ds^{-1} := (1+D^2)^{-\frac{1}{2}}$ is a compact operator.

In Connes' quantised calculus the compact operator $\ds^{-1}$
is the analogue of an inf\/initesimal length element $ds$~\cite[p.~545]{CN}, \cite[p.~157]{C4},
and the integral of $a \in \A$ is represented by the operator-theoretic expression~\cite[p.~545]{CN}, \cite[p.~158]{C4},
\begin{equation} \label{eq:ncg_int}
\Tr_\omega(a\ds^{-n})  , \qquad a \in \A,
\end{equation}
where $n$ (assumed to exist) is the smallest value such that $\ds^{-p}$ is trace class for $p > n$, and the positive trace $\Tr_\omega$
(nominally a `Dixmier trace')
measures the log divergence of the trace of~$a\ds^{-n}$.  The link between integration in the traditional sense and sums of eigenvalues of the product $a\ds^{-n}$ is not obvious.

Segal, in \cite{Se2}, formalised an operator-theoretic view of integration theory.  The subsequently developed theory\footnote{The explanation here is for semif\/inite von Neumann algebras with separable pre-dual.} \cite{Kunz1958, Stinespring1959}, through to \cite{Nelson1974, FKsnum,PX2003}, features the notion of a faithful normal trace $\tau$ on
a~$*$-algebra $\mathcal{N}$ of bounded linear operators on $H$ that is closed in the weak operator topology (a von Neumann algebra).
The trace $\tau$ is used to build an $L^p$-theory with $L^\infty(\mathcal{N},\tau) = \mathcal{N}$
and $L^p(\mathcal{N},\tau)$ being the completion of the set
$\{ x \in \mathcal{N} | \tau(|x|^p) < \infty \}$.  The trace (and the integration theory) is called f\/inite
if $\tau(1_\mathcal{N}) < \infty$ where $1_\mathcal{N}$ is the identity of $\mathcal{N}$.
The integration theory associated to the pair $(\mathcal{N},\tau)$ has
analogous machinery to measure theory, including
Radon--Nikodym theorems, dominated and monotone convergence theorems\footnote{The condition of normality of a linear functional preempts the monotone convergence theorem.  If $\mathcal{N}$
is a closed $*$-subalgebra of $\mathcal{L}(H)$ in the weak operator topology then
a linear functional $\rho \in \mathcal{N}^*$ is normal if
$\rho(S) = \sup_\alpha \rho(S_\alpha)$ for all
increasing nets of positive operators $\{ S_\alpha \} \subset \mathcal{N}$ in the strong operator topology with l.u.b.~$S$.  This condition is equivalent
to $\rho(\sum_i P_i) = \sum_i \rho(P_i)$ for all sets $\{ P_i \}$ of pairwise orthogonal projections
from $\mathcal{N}$~\cite[p.~67]{Ped}, i.e.~analogous to additivity of a measure on a
$\sigma$-algebra of sets.}, etc.
Without the condition of normality this machinery fails.

{\sloppy A commutative example of the theory (and the prototype) is a regular
Borel measure space~$(X,\mu)$.  For the von Neumann
algebra
of complex functions $f \in L^\infty(X,\mu)$ acting on
the separable Hilbert space $L^2(X,\mu)$ by (almost everywhere)
pointwise multiplication, the integral
\[
\tau_\mu(f) = \int_X f d\mu  , \qquad f \in L^\infty(X,\mu)
\]
is a faithful normal trace.  The noncommutative $L^p$-spaces are
the usual $L^p$-spaces, i.e.
\[
L^p(L^\infty(X,\mu),\tau_\mu) \equiv L^p(X,\mu).
\]

}

A noncommutative example is
the von Neumann algebra $\mathcal{L}(H)$ of bounded
linear operators on a separable Hilbert space~$H$.  The canonical trace
\[
\Tr(T) = \sum_{m=1}^\infty  \inprod{h_m}{Th_m}   , \qquad T \in \mathcal{L}(H),
\]
where $\{ h_m \}_{m=1}^\infty$ is an orthonormal basis
of $H$, is a faithful normal trace.
The resulting $L^p$-spaces
\begin{gather} \label{eq:schatten}
\mathcal{L}^p(H) := \{ T \in \mathcal{L}(H) | \Tr(|T|^p) < \infty \}
\end{gather}
are the Schatten--von Neumann ideals of compact operators \cite{S} (\cite{JvN1937} for f\/inite dimensional~$H$).
A consequence of the Radon--Nikodym theorem
in this context is that for any normal linear functional $\phi$
on a von Neumann algebra $\mathcal{N} \subset \mathcal{L}(H)$
there is a trace class `density' $T \in \mathcal{L}^1(H)$ such that, \cite[p.~55]{Ped},
\begin{equation} \label{eq:ncg_int_standard}
\phi(a) = \Tr(aT) , \qquad \fa \, a \in \mathcal{N} .
\end{equation}

Combining the examples, for every regular Borel measure $\mu$
on a space $X$ there is a positive trace class operator $T_\mu$ on the Hilbert space $L^2(X,\mu)$ such that
\begin{equation} \label{eq:ncg_int_example}
\int_X f d\mu = \Tr( f T_\mu) , \qquad f \in L^\infty(X,\mu) .
\end{equation}
This formula exhibits the standard association between integration
and the eigenvalues of the product $f T_\mu$, a compact linear operator.

What are the dif\/ferences between \eqref{eq:ncg_int}, the integral in NCG, and \eqref{eq:ncg_int_standard}, the integral according to standard noncommutative integration theory?
Let us mention basic facts about \eqref{eq:ncg_int}.

Firstly, since $\ds^{-n}$ belongs to the ideal of compact linear operators
\begin{equation} \label{eq:intro_dix_ideal}
\mathcal{M}_{1,\infty}(H) := \left\{ T \in \mathcal{L}^\infty(H) \, \Big| \,
\sup_{n \in \NN} \frac{1}{\log(1+n)} \sum_{k=1}^n \mu_k(T)  < \infty \right\}
\end{equation}
(where $\mathcal{L}^\infty(H)$ denotes the compact linear operators and $\mu_{n}(T)$, $n \in \NN$, are the eigenvalues of $|T|=(T^*T)^{1/2}$ arranged in decreasing order),
and the positive traces $\Tr_\omega$ are linear functionals on $\mathcal{M}_{1,\infty}$, the formula
\begin{equation} \label{eq:ncg_ext}
\Tr_\omega(S\ds^{-n})   , \qquad S \in \mathcal{L}(H)
\end{equation}
def\/ines a linear functional on $\mathcal{L}(H)$.
In particular, by restriction \eqref{eq:ncg_ext} is a linear functional on the $*$-algebra~$\A$,
on the closure~$A$ of $\A$ in the uniform operator topology (a separable $C^*$-subalgebra of $\mathcal{L}(H)$ since $\A$ is countably generated), and on
the closure $\A''$ of $\A$ in the weak (or strong) operator topology (the bicommutant, a von Neumann subalgebra of $\mathcal{L}(H)$ \cite[\S~2.4]{BR}).

Secondly, if $0 < S \in \mathcal{L}(H)$,
\[
 \Tr_\omega(S\ds^{-n}) =  \Tr_\omega\big(\sqrt{S}\ds^{-n}\sqrt{S}\big) \geq 0 .
\]
Hence \eqref{eq:ncg_ext} def\/ines a positive linear functional
on $\A$, $A$ or $\A''$.

Thirdly,
\[
| \Tr_\omega(S\ds^{-n}) | \leq C \nm{S}   , \qquad \fa \, S \in \mathcal{L}(H)
\]
for a constant $C>0$.  Hence the same formula~\eqref{eq:ncg_ext} is the unique uniformly continuous extension to~$A$ and
$\A''$ of the positive linear functional~\eqref{eq:ncg_int} on $\A$.  Denote this linear functional by $0 < \Phi_\omega \in (\A'')^*$.

Fourthly, $\Phi_\omega$ is a trace on $\A''$  \cite[Theorem~1.3]{Cip1996},
\cite[p.~313]{CN}, \cite[p.~280]{CRSS}, \cite[p.~47]{FB}.

The fundamental dif\/ference between the functionals \eqref{eq:ncg_int_standard} and \eqref{eq:ncg_ext} is that $\Tr_\omega(S\ds^{-n}) = 0$ for all positive f\/inite rank operators $S$ (see Section~\ref{sec:dtr} below).  In comparison $\Tr(ST) > 0$ for a~trace class operator $0 < T \in \mathcal{L}^1(H)$ with trivial kernel and every f\/inite rank operator $S>0$.
This property of vanishing on the f\/inite rank operators, called singularity, is essential for the fundamental theorem of integration theory in
NCG (Connes' trace theorem \cite{C3}, and Theorem~\ref{thm:2}
in this text).  The property of singularity, however, implies that
the functional~\eqref{eq:ncg_ext} is either trivial or else is a non-normal functional on
$\mathcal{L}(H)$. To see this, let $S_n$ be a sequence
of positive f\/inite rank operators converging upward strongly
to the identity~$1$ of $\mathcal{L}(H)$ and suppose
\eqref{eq:ncg_ext} is a normal linear functional on $\mathcal{L}(H)$, then
\[
0 = \sup_n \Tr_\omega(S_n\ds^{-n}) =  \Tr_\omega(\ds^{-n}).
\]
Hence either $\Tr_\omega(\ds^{-n})=0$ or else~\eqref{eq:ncg_ext} is non-normal.

The property of singularity does not exclude the trace $0 < \Phi_\omega \in (\A'')^*$ (the restriction of~\eqref{eq:ncg_ext} to the
von Neumann algebra $\A'' \subset \mathcal{L}(H)$) from being a normal linear functional on~$\A''$.  There are spectral triples where $\Phi_\omega$
is not a normal linear functional on~$\A''$
(evidently when~$\A''$ contains a f\/inite rank operator), and where $\Phi_\omega$ is a normal linear functional on~$\A''$
(see Theorem~\ref{thm:result2}).

Currently we lack, in general, a characterisation of the relationship between $D$
and $\A$ that implies $\Phi_\omega \in (\A'')^*$ is a normal linear functional.  All we know of is a suf\/f\/icient condition, due to the authors, see Theorem~\ref{thm:normM} in Section~\ref{sec:char:dct}.  It has been suggested that $\A''$ containing no f\/inite rank operators is
a necessary and suf\/f\/icient condition for normality of $\Phi_\omega$.  This conjecture is open.  Non-normality of the linear functional on $\mathcal{L}(H)$ def\/ined by \eqref{eq:ncg_ext} can have advantages, as indicated by Connes \cite[p.~326]{CN}.
However, normality of the trace $0 < \Phi_\omega \in (\A'')^*$ is coincident
with our ability to apply the standard tools of integration theory to the formula~\eqref{eq:ncg_int}.

In general, the value of the trace $0 < \Phi_\omega \in (\A'')^*$ on the projections of $\A''$,
\begin{equation} \label{eq:ncg_meas}
\Phi_\omega(E)   , \qquad E^* = E = E^2 \in \A'' ,
\end{equation}
is only analogous to an element of the ba space of f\/inite and f\/initely additive functionals on a~Borel $\sigma$-algebra \cite[IV.2]{DSI}.
The integral~\eqref{eq:ncg_int} is constructed
by limits of f\/inite linear combinations of the values in~\eqref{eq:ncg_meas} if and only if the trace $\Phi_\omega$ is normal.

In this review, we discuss the known ways to measure the `log divergence of the trace', including residue formulas, heat kernels and Dixmier traces.  Each of these ways may result
in dif\/ferent functionals in~\eqref{eq:ncg_int}.
It is desirable to know when all these ways agree, and that \eqref{eq:ncg_int}
def\/ines a unique linear functional on~$\A$.
There are a variety of dif\/f\/iculties behind this problem.  For example, if $\mathcal{V}$
denotes a set of linear functionals on $\mathcal{M}_{1,\infty}$ `measuring the log divergence of the trace', the set
\[
\mathcal{K}_\mathcal{V} := \{ T \in \mathcal{M}_{1,\infty} \, |\,
\mathrm{f}(T) = \mathrm{const} \  \fa \, \mathrm{f} \in \mathcal{V} \}
\]
is not an ideal \cite[Remark~3.6]{AF}.  Therefore, $\ds^{-n} \in \mathcal{K}_\mathcal{V}$
will not imply $a \ds^{-n} \in \mathcal{K}_\mathcal{V}$ for all $a \in \A$ in general.  Specifying
$\ds^{-1}$ alone is not suf\/f\/icient for a unique value from the formula~\eqref{eq:ncg_int}.
Additionally, if one picks smaller sets, $\mathcal{V}$, or larger sets, $\mathcal{V}'$, of linear functionals, it may be that $\mathcal{K}_{\mathcal{V}} \not= \mathcal{K}_{\mathcal{V}'}$ (discussed in Section~\ref{sec:meas}).  Allied to this problem is the fact that, if $a \in \A$ is a selfadjoint operator and $a\ds^{-n} \in \mathcal{K}_\mathcal{V}$ it is not necessarily true that the decompositions $a_+ \ds^{-n}, a_- \ds^{-n} \in \mathcal{K}_\mathcal{V}$, where $a_+$ and $a_-$ are the positive and negative parts of~$a$.

Currently, in general, there
is no characterisation, given $\ds^{-n}$, what relationship with
\mbox{a~$*$-al\-gebra} $\A$ is necessary so that $\Phi_\omega(a) = \Tr_\omega(a\ds^{-n})$ specif\/ies a unique trace (meaning the value of $\Phi_\omega(a)$ is independent of the linear functional $\Tr_\omega \in \mathcal{V}$).
This is an open problem.  Checking uniqueness is a case by case basis which requires identifying $\Tr_\omega(a\ds^{-n})$ with a~known functional, e.g.~the Lebesgue integral for a compact Riemannian manifold \cite{C3,FB}, or the Hausdorf\/f measure of Fuchsian circles
and self-similar fractals \cite[IV.3]{CN}, \cite{GI2003}, or estimating the singular values of the product $a\ds^{-n}$.

In the review we shall provide some suf\/f\/icient
conditions for normality of the trace $0 < \Phi_\omega \in (\A'')^*$.
We review the demonstration of normality for the case of a compact Riemannian manifold (in Section~\ref{sec:origin}) and the noncommutative torus (in Section~\ref{sec:char}).

Finally, we discuss the role of other positive singular traces on $\mathcal{M}_{1,\infty}$, besides Dixmier traces,
in the formula~\eqref{eq:ncg_int}.  We show that, even in Connes' trace theorem \cite[Theorem~1]{C3}, there exist other positive singular traces $\rho$ on the compact operators, which are not Dixmier traces, such that $\rho(P) = \mathrm{Res}(P)$ (up to a positive constant) for a classical pseudo-dif\/ferential $P$ of order $-n$ on a $n$-dimensional compact Riemannian manifold, where~$\mathrm{Res}$ is the Wodzicki residue.

\subsection*{Notation}

{\sloppy
Throughout $H$ denotes a separable complex Hilbert space,
$\mathcal{L}(H)$ (resp.~$\mathcal{L}^\infty \equiv \mathcal{L}^\infty(H)$, $\mathcal{L}^p \equiv \mathcal{L}^p(H)$, $p \in [1,\infty)$)
denotes the linear bounded (resp.~compact, Schatten--von Neumann (see~\eqref{eq:schatten})) operators on~$H$.  Throughout the text
$\mathcal{M}_{1,\infty} \equiv \mathcal{M}_{1,\infty}(H)$
denotes the ideal of compact operators def\/ined in~\eqref{eq:intro_dix_ideal} or~\eqref{eq:dixmier_ideal} below. Denote by $\ell^\infty$ (resp.~$c$, $c_0$, $\ell^p$)
the bounded (resp.~convergent, convergent to zero, $p$-summable) sequences.  Similarly, for a regular Borel measure space
$(X,\mu)$, $L^\infty(X,\mu)$ (resp.~$L^p(X,\mu)$,
$C_b(X)$, $C_0(X)$) denotes the essentially bounded (resp.~$p$-integrable, continuous and bounded, continuous and vanishing at inf\/inity) (equivalence classes of) functions on~$X$.
When $X$ is a Riemannian manifold and $dx$ is the volume form on~$X$, we write
$L^p(X) \equiv L^p(X,dx)$.

}

\section{Dixmier traces} \label{sec:dtr}

For a separable complex Hilbert space $H$, denote by $\mu_{n}(T)$, $n \in \NN$, the singular values of a~compact operator $T$ (the singular values are the eigenvalues of the operator
$|T| = (T^*T)^{1/2}$ arranged with multiplicity in decreasing order  \cite[\S~1]{S}).  Def\/ine the logarithmic average
\begin{equation} \label{eq:alpha_first}
\alpha( \{ \mu_n(T) \}_{n=1}^\infty ) := \left\{ \frac{1}{\log(1+k)} \sum_{n=1}^k \mu_n(T) \right\}_{k=1}^\infty , \qquad T \in \mathcal{L}^\infty.
\end{equation}
Then
\begin{equation} \label{eq:dixmier_ideal}
\mathcal{M}_{1,\infty}:= \mathcal{M}_{1,\infty}(H) = \inset{T}{\nm{T}_{1,\infty}:=\sup_k \alpha( \{ \mu_n(T) \}_{n=1}^\infty )_k < \infty}
\end{equation}
def\/ines a Banach ideal of compact operators whose sequence of partial sums of singular values is of order logarithm.
We refer to the recent paper of Pietsch~\cite{Pie}, discussing the origin of this object in mathematics.
In \cite{Dix}, J.~Dixmier constructed a non-normal semif\/inite trace  (a~Dixmier trace) on
$\mathcal{L}(H)$ using the weight
\begin{equation} \label{eq:dtr0}
\Tr_\omega(T) := \omega \left( \left\{ \frac{1}{\log(1+k)} \sum_{n=1}^k \mu_n(T)\right\}_{k=1}^\infty  \right)  , \qquad T>0,
\end{equation}
where $\omega$ is a state on $\ell^\infty$ associated to a translation and dilation invariant state
on $\RR$.
The sequence \eqref{eq:alpha_first},
while bounded for $T \in \mathcal{M}_{1,\infty}$,
is not generally convergent.  There are operators $0 < T \in \mathcal{M}_{1,\infty}$ for which the sequence $\alpha( \{ \mu_n(T) \}_{n=1}^\infty ) \in c$ (called Tauberian operators),
and for the Tauberian operators we want a Dixmier trace to be equal to the limit
$\lim\limits_{k \to \infty} \alpha( \{ \mu_n(T) \}_{n=1}^\infty )_k$.
However,
the linear span of the set of Tauberian operators does not form an
ideal of $\mathcal{M}_{1,\infty}$.  For the semif\/inite domain of
a Dixmier trace to be an ideal of compact operators, the use
of a~state~$\omega$ on~$\ell^\infty$ such that
$\omega(a_k) = \lim\limits_k a_k$ when $a_k \in c$
(known as a generalised limit) cannot be avoided.
Note that $|\omega(a_k)| \leq \nm{a_k}_\infty$,
hence
\begin{equation} \label{eq:dtr_cont}
| \Tr_\omega(T) | \leq \nm{T}_{1,\infty} , \qquad T \in \mathcal{M}_{1,\infty} .
\end{equation}

A Dixmier trace $\Tr_\omega$ is positive by construction
and continuous by \eqref{eq:dtr_cont}, i.e.~$0 < \Tr_\omega \in \mathcal{M}_{1,\infty}^*$.  That it is
a trace follows from the identity $\mu_n(U^*TU) = \mu_n(T)$ for any unitary $U \in \mathcal{L}(H)$.  That it is singular
(vanishes on f\/inite rank operators, in fact on any trace class operator) follows from
$\alpha( \{ \mu_n(S) \}_{n=1}^\infty ) \in c_0$ for any f\/inite rank operator $S$ (resp.~$S \in \mathcal{L}^1$).
The basic facts stated in the introduction now follow (for the fourth see the cited references).
For the third notice that $\nm{ST}_{1,\infty} \leq \nm{S} \nm{T}_{1,\infty}$ (hence $|\Tr_\omega(ST)| \leq \nm{S} \nm{T}_{1,\infty}$) for all $S \in \mathcal{L}(H)$, $T \in \mathcal{M}_{1,\infty}$.

The non-trivial feature of $\Tr_\omega$ is linearity.  Dixmier's specif\/ication that $\omega$
is constructed from a translation and dilation invariant state on $\RR$ provides suf\/f\/icient conditions for linearity.  Later formulations of the Dixmier trace considered wider or narrower specif\/ications suited to their context.

Unfortunately, the specif\/ications are a mess.  We will now introduce a wide range of sets of linear functionals, all of which
play some part in known results concerning the Dixmier trace.

Let $S(\ell^\infty)$ denote the states on $\ell^\infty$, i.e.~the positive linear functionals $\sigma$ such that $\sigma((1,1,\dots))$ $ =\, 1$.  Let $S([1,\infty))$ (resp.~$S([0,\infty))$) denote the states of $L^\infty([1,\infty))$ (resp.~$L^\infty([0,\infty))$).  A state on $\ell^\infty$ (resp.~$L^\infty([1,\infty))$, $L^\infty([0,\infty))$) is singular if it vanishes on f\/inite sequences (resp.~functions with a.e.~compact support).  It is easy to show a state $\sigma$ on
$\ell^\infty$ (resp.~$L^\infty([1,\infty))$, $L^\infty([0,\infty))$)
is singular if and only if
\begin{gather*}
\liminf_{n \to \infty} a_n \leq \sigma( \{a_n \}_{n=1}^\infty) \leq \limsup_{n \to \infty} a_n   , \qquad a_n > 0 \quad \forall\, n \in \NN
\\
\mbox{(resp.} \ \
\mathrm{ess\hbox{-}}\liminf_{t \to \infty} f(t) \leq \sigma( f ) \leq \mathrm{ess\hbox{-}}\limsup_{t \to \infty} f(t) \ , \quad f > 0).
\end{gather*}
That is, a state is singular if and only if it is an extension to
$\ell^\infty$ (resp.~$L^\infty([1,\infty))$, $L^\infty([0,\infty))$) of the ordinary limit of convergence sequences (resp.~limit at inf\/inity of convergent at inf\/inity functions).  The Hahn--Banach theorem guarantees that these linear extensions, called generalised limits, exist.  Denote by $S_\infty(\ell^\infty)$ (resp.\ $S_\infty([1,\infty))$, $S_\infty([0,\infty))$) the set of singular states (or genera\-lised limits).  Denote
by $SC_\infty([1,\infty))$ (resp.\ $SC_\infty([0,\infty))$) the same construct with
functions from $L^\infty([1,\infty))$ (resp.\ $L^\infty([0,\infty))$) replaced by continuous
functions from $C_b([1,\infty))$ (resp.\ $C_b([0,\infty))$).

Let $\ceil{x}$, $x \geq 0$, denote the ceiling function.
For $j \in \NN$ def\/ine the maps $\ell^\infty \to \ell^\infty$ by
\begin{gather}
T_j( \{a_k\}_{k=1}^\infty )   :=   \{ a_{k+j} \}_{k=1}^\infty, \label{eq:def_Td} \\
D_j( \{a_k \}_{k=1}^\infty)   :=   \{ a_{\ceil{j^{-1}k}} \}_{k=1}^\infty, \label{eq:def_Dd} \\
C( \{a_k \}_{k=1}^\infty)    :=
\left\{ \frac{1}{n} \sum_{k=1}^n  a_k \right\}_{n=1}^\infty \label{eq:def_C} .
\end{gather}
Def\/ine subsets of $S_\infty(\ell^\infty)$ by
\begin{gather}
BL   :=   \{ \omega \in S(\ell^\infty) \, | \,
\omega \circ T_j = \omega \ \fa\,  j \in \NN\}, \label{eq:def_BL} \\
CBL   :=   \{ \omega \in S(\ell^\infty) \, | \,
\omega = \sigma \circ C, \ \sigma \in S_\infty(\ell^\infty) \},
\label{eq:def_CBL} \\
B(C)   :=   \{ \omega \in S(\ell^\infty) \, |\,
\omega \circ C = \omega \}, \label{eq:def_B(C)} \\
D_2   : =   \{ \omega \in S_\infty(\ell^\infty) \, | \,
(\omega \circ D_2 - \omega)(\alpha( \{ \mu_n(T) \}_{n=1}^\infty )) = 0 \ \forall\,  T\in\mathcal{M}_{1,\infty} \}, \label{eq:def_D2} \\
DL   :=   \{ \omega \in S(\ell^\infty) \, |\,
\omega \circ D_j = \omega \ \fa \, j \in \NN\}, \label{eq:def_DL}
\end{gather}
where $\alpha$ is from \eqref{eq:alpha_first}. The set $BL$ is the classical set of Banach limits \cite{li-2,li-3}.
The elements of $CBL$ are called Ces\`{a}ro--Banach limits and the
elements of $B(C)$ are called Ces\`{a}ro invariant Banach limits.
We note that (as the names suggest)
\begin{equation} \label{eq:disc_inclusions}
B(C) \subsetneq CBL \subsetneq BL
\end{equation}
and
\begin{equation} \label{eq:disc_inclusions2}
DL \subsetneq D_2 .
\end{equation}

For $a > 0$ def\/ine the maps $L^\infty([0,\infty)) \to L^\infty([0,\infty))$ by
\begin{gather}
T_a(f)(t)   :=   f(t+a), \label{eq:def_Tc} \\
D_a(f)(t)   :=   f(a^{-1}t), \label{eq:def_Dc} \\
P_a(f)(t)   :=   f(t^{a}), \label{eq:def_Pc} \\
C(f)(t)   :=   \frac{1}{t} \int_0^t f(s)ds. \label{eq:def_Cc}
\end{gather}
Def\/ine $L^{-1} : L^\infty([1,\infty)) \to L^\infty([0,\infty))$ by
\[
L^{-1}(g)(t) := g(e^t)
\]
and, if $G : L^\infty([0,\infty)) \to L^\infty([0,\infty))$, set
\begin{equation} \label{eq:def_L_maps}
L(G) : \ \ L^\infty([1,\infty)) \to L^\infty([1,\infty))   , \qquad G \mapsto L \circ G \circ L^{-1} .
\end{equation}
Similarly, if $\sigma \in S([0,\infty))$ set
\begin{equation} \label{eq:def_L_states}
L(\sigma) \in S([1,\infty))   , \qquad \sigma \mapsto \sigma \circ L^{-1} .
\end{equation}
These def\/initions imply $L(\sigma) \circ L(G) = L(\sigma \circ G)$.
It is known $L(T_{a}) = D_{e^{-a}}$ and $L(D_{a}) = P_{a^{-1}}$, see  \cite[\S~1.1]{CPS}.  Def\/ine subsets of $S_\infty([0,\infty))$ and $S_\infty([1,\infty))$ by
\begin{gather}
BL[0,\infty)   :=   \{ \omega \in S([0,\infty)) \, |\,
\omega \circ T_a = \omega \ \fa \, a > 0\}, \label{eq:def_BLi} \\
CBL[0,\infty)   :=   \{ \omega \in S([0,\infty))\, |\,
\omega = \sigma \circ C, \ \sigma \in S_\infty([0,\infty)) \},
\label{eq:def_CBLi} \\
B(C)[0,\infty)   :=   \{ \omega \in S([0,\infty)) \, |\,
\omega \circ C = \omega \}, \label{eq:def_B(C)i} \\
DL[0,\infty)   :=   \{ \omega \in S([0,\infty)) \, |\,
\omega \circ D_a = \omega \ \fa \, a > 0\}, \label{eq:def_DL0i} \\
DL[1,\infty)   :=   \{ L(\omega) \, |\,
\omega \in BL[0,\infty) \}, \label{eq:def_DLi} \\
CDL[1,\infty)   :=   \{ L(\omega) \, |\,
\omega \in CBL([0,\infty)) \},
\label{eq:def_CLi} \\
D(C)[1,\infty)   :=   \{ L(\omega) \, |\,
\omega \in B(C)([0,\infty)) \}, \label{eq:def_D(C)i} \\
PDL[1,\infty)   :=   \{ L(\omega) \, |\,
\omega \in BL[0,\infty) \cap DL[0,\infty) \} ,\label{eq:def_PDLi} \\
MDL[1,\infty)    :=   \{ L(\omega) \, |\,
\omega \in BL[0,\infty) \cap DL[0,\infty) \cap B(C)[0,\infty)  \} . \label{eq:def_MDLi}
\end{gather}
A subscript of $c$ for any of the sets \eqref{eq:def_BLi}--\eqref{eq:def_MDLi} denotes the replacement of
$L^\infty$-functions by $C_b$-functions, i.e.~$BL_c[0,\infty)$ denotes
the translation invariant states on the $C^*$-algebra $C_b([0,\infty))$
and it is a subset of $SC_\infty([0,\infty))$.

We now def\/ine the Dixmier traces corresponding to the singular states
\eqref{eq:def_BL}--\eqref{eq:def_D2}, \eqref{eq:def_BLi}--\eqref{eq:def_MDLi} and subscript~$c$ versions.

Let $g^*$ be the monotonically decreasing left-continuous rearrangement of $|g|$, where $g \in L^\infty([0,\infty))$~\cite{KPS1982}. Def\/ine
\begin{equation} \label{eq:alpha}
\alpha_t(g^*) := \frac{1}{\log(1+t)} \int_1^t g^*(s)ds  , \qquad t \geq 1
\end{equation}
(the continuous version of \eqref{eq:alpha_first}).
Def\/ine the f\/loor map $p$ from~$\ell^\infty$ to bounded Borel functions on $[0,\infty)$ (or $[1,\infty)$) by
\begin{equation} \label{eq:def_p}
p( \{a_k \}_{k=1}^\infty)(t) :=
\sum_{k=1}^\infty  a_k \chi_{[k,k+1)}(t)   , \qquad t \geq 0 \ \ \mbox{(or $1$)}
\end{equation}
and the restriction map $r$ from bounded Borel functions on $[0,\infty)$ (or $[1,\infty)$) to $\ell^\infty$
\begin{equation} \label{eq:def_r}
r(f)_n := f(n)    , \qquad n \in \mathbb{N} .
\end{equation}
We have the following sets of unitarily invariant weights
on $\mathcal{M}_{1,\infty}$:
\begin{gather}
\mathcal{V}_B  =  \{ \tr_\omega(T):=
\omega \circ \alpha_t \circ p( \{ \mu_n(T) \} ), T>0  \, |\, \omega \in B\}, \nonumber\\
\phantom{\mathcal{V}_B  =}{} \ B = \left\{  \begin{array}{l}
DL[1,\infty), \ DL_c[1,\infty), \\
CDL[1,\infty), \  CDL_c[1,\infty), \\
D(C)[1,\infty), \ D(C)_c[1,\infty), \\
PDL[1,\infty), MDL[1,\infty)
\end{array} \right. \label{eq:def_dt_1} \\
\mathcal{V}_B  =  \{ \Tr_\omega(T):=
\omega \circ r \circ \alpha_t \circ p( \{ \mu_n(T) \}), T>0  \, |\, \omega \in B\}, \qquad B =
D_2, \label{eq:def_dt_2} \\
\mathcal{W}_Y  =  \{ \mathrm{f}_\omega(T):=
\omega \circ \alpha_{e^t} \circ p( \{ \mu_n(T) \}), T>0 \, |\, \omega \in Y\}, \nonumber\\
\phantom{\mathcal{W}_Y  =}{} \ Y = \left\{ \begin{array}{l}
BL[0,\infty), \ BL_c[0,\infty), \\
CBL[0,\infty), \  CBL_c[0,\infty), \\
B(C)[0,\infty), \ B(C)_c[0,\infty)
\end{array} \right. \label{eq:def_dt_3} \\
\mathcal{W}_Y  =  \{ \mathrm{F}_\omega(T):=
\omega \circ r \circ \alpha_{e^t} \circ p( \{ \mu_n(T) \}), T>0 \, |\, \omega \in Y\}, \!\!\qquad  Y =
BL, CBL, B(C) .\!\!\! \label{eq:def_dt_4}
\end{gather}
By linear extension each of the unitarily invariant weights def\/ines a positive trace with
do\-main~$\mathcal{M}_{1,\infty}$.

The reader should rightly be bewildered by this variety of Dixmier traces.  Before we explain the origin of so many options,
let us clarify the picture.  From Appendix~\ref{sec:app_A} we have
\begin{equation} \label{fig:simplify}
\begin{array}{ccccc}
\mathcal{V}_{D_2} & \supset &
\begin{array}{ccc}
%\mathcal{V}_{DL} =
\ \ \ \ \ \ \ \ \ \ \ \  \mathcal{V}_{DL[1,\infty)}
= \mathcal{V}_{DL_c[1,\infty)} \\
\shortparallel \\
\mathcal{W}_{BL} = \mathcal{W}_{BL[1,\infty)}
= \mathcal{W}_{BL_c[1,\infty)}
\end{array}				 & \supset & \mathcal{V}_{PDL_{[1,\infty)}} \\[24pt]
& & \cup & & \\[4pt]
& & \begin{array}{ccc}
\ \ \ \ \ \ \ \ \ \ \ \  \mathcal{V}_{CDL[1,\infty)}
= \mathcal{V}_{CDL_c[1,\infty)} \\
\shortparallel \\
\mathcal{W}_{CBL} = \mathcal{W}_{CBL[1,\infty)}
= \mathcal{W}_{CBL_c[1,\infty)}
\end{array} & & \\[24pt]
& & \cup & & \\[4pt]
& & \begin{array}{ccc}
\ \ \ \ \ \ \ \ \ \ \ \ \mathcal{V}_{D(C)[1,\infty)}
= \mathcal{V}_{D(C)_c[1,\infty)} \\
\shortparallel \\
\mathcal{W}_{B(C)} = \mathcal{W}_{B(C)[1,\infty)}
= \mathcal{W}_{B(C)_c[1,\infty)}
\end{array} &  \supset & \mathcal{V}_{MDL_{[1,\infty)}}
\end{array}
\end{equation}
Why did such a variety result?  Connes' trace theorem~\cite{C3} used the generalised limits origi\-nally specif\/ied by Dixmier~\cite{Dix}.
Connes later introduced Dixmier traces def\/ined by the sets $CDL_c[1,\infty)$ \cite[pp.~303--308]{CN},
after explaining~\cite[p.~305]{CN}, following Dixmier~\cite{Dix}, that any $\omega \in D_2$
was suf\/f\/icient to def\/ine a Dixmier trace.  In the same location~\cite[p.~308]{CN} Connes introduced the idea of measurability associated to the set of traces def\/ined by $CDL_c[1,\infty)$ (see the next section).

Translation invariance in Dixmier's original specif\/ication is redundant, hence the set of functionals $\mathcal{V}_{DL[1,\infty)}$ are generally termed the (original) Dixmier traces, while the more restrictive set $\mathcal{V}_{CDL_c[1,\infty)}$ we termed Connes--Dixmier traces, in~\cite{LSS}, to distinguish them. The text~\cite{LSS} also initiated the identif\/ications in \eqref{fig:simplify}, where our aim was to identify Dixmier and Connes--Dixmier traces with subsets of singular symmetric functionals \cite{DPSS,DPSSS1,DPSSS2}, i.e.~$\mathcal{V}_{DL[1,\infty)}=\mathcal{W}_{BL}$
and $\mathcal{V}_{CDL_c[1,\infty)} = \mathcal{W}_{CBL}$.
The form of the functionals
$\mathcal{W}_{BL, CBL}$ is simpler than it looks in~\eqref{eq:def_dt_4}.  Indeed, as a result of~\cite{LSS} and~\cite{DPSSS1}, one can write any Dixmier trace (up to a constant) as
\[
\mathrm{F}_\omega(T) = \omega \left( \left\{ \frac{1}{N} \sum_{k=1}^{2^N} \mu_k(T) \right\}_{N=1}^\infty \right)   , \qquad T > 0
\]
for some Banach limit $\omega$, and any Connes--Dixmier trace
using the same formula except with a~Ces\`{a}ro--Banach limit $\omega \in CBL$.

The sets $\mathcal{V}_{D(C)[1,\infty)}$ and $\mathcal{V}_{MDL[1,\infty)}$
were introduced with a dif\/ferent intention.  The paper~\cite{CPS}
by A.~Carey, J.~Phillips, and the second author,
demonstrated the residue and heat kernel
formulations of the Dixmier trace, and did so with the type I von Neumann algebra $(\mathcal{L}(H),\Tr)$
replaced with any semif\/inite von Neumann algebra $(\mathcal{N},\tau)$
with faithful normal semif\/inite trace~$\tau$, see also \cite{AF, GI1995,FB}.
The generalised $s$-numbers $\rho_s(T)$ of a $\tau$-compact operator $T$ in $\mathcal{N}$,~\cite{FKsnum}, provide a
generalisation of singular values, such that $\rho_s(T) = p(\{\mu_n(T)\}_{m=1}^\infty)$ for the pair $(\mathcal{L}(H),\Tr)$.
To use the generalised $s$-numbers, one replaces the sequence \eqref{eq:alpha_first} by the function \eqref{eq:alpha}.  For the
semif\/inite pair $(\mathcal{N},\tau)$ only the functionals $\mathcal{V}_{DL[1,\infty)}$ (and all continuous variants in~\eqref{eq:def_dt_1} and~\eqref{eq:def_dt_3}) can be def\/ined.

The set $MDL[1,\infty)$, which is non-empty \cite[p.~74]{CPS}, was a constraint used in \cite{CPS} to identify the Dixmier trace with residues of a zeta function and heat kernel asymptotics
(see Section~\ref{sec:res} below).  Later results on zeta functions~\cite{CRSS}  enabled a residue formulation for all Dixmier traces in the set $\mathcal{V}_{PDL[1,\infty)}$.  It is still open whether a residue formulation exists for
all functionals in the set $\mathcal{V}_{DL[1,\infty)}$.
While residue formulations were able to drop the Ces\`{a}ro invariance condition (which is unnecessary for the weak$^*$-Karamata theorem~\cite[p.~271]{CRSS}), a later heat kernel formulation due to A.~Sedaev dropped the power invariance requirement (as  \cite[Proposition~4.3]{CRSS} is unnecessary for the heat kernel).
In~\cite{Sed2009}, heat kernel formulations were derived
for Connes--Dixmier traces in the sets $\mathcal{V}_{CDL_c[1,\infty)}$ and $\mathcal{V}_{D(C)_c[1,\infty)}$.
%Due to the nature of the heat kernel transform the reference~\cite{Sed2009} probably represents the f\/inal form of heat kernel formulations of the Dixmier trace, see also the recent
%pre-print concerning Lidskii formulas~\cite{ssz2010}.

The set $D(C)[1,\infty)$ was used in~\cite[\S~1.3]{FB}. %(see~\cite[\S~1.3]{FB}).
Finally, the set $B(C)$ was introduced in the paper~\cite{SemSuk2009} with an application to measurability (see the next section).

It will be convenient to simplify the schematic \eqref{fig:simplify} of Dixmier traces.  Set
\begin{gather}
\mathcal{V}_{\frac{1}{2}}   :=   \mathcal{V}_{D_2}, \label{eq:dtr0.5} \\
\mathcal{V}_{1}   :=
\mathcal{V}_{DL[1,\infty)}
= \mathcal{V}_{DL_c[1,\infty)} = \mathcal{W}_{BL_c[1,\infty)} = \mathcal{W}_{BL[1,\infty)} = \mathcal{W}_{BL},
\label{eq:dtr1} \\
\mathcal{V}_{2}   :=   \mathcal{V}_{CDL[1,\infty)}
= \mathcal{V}_{CDL_c[1,\infty)} = \mathcal{W}_{CBL_c[1,\infty)} = \mathcal{W}_{CBL[1,\infty)} = \mathcal{W}_{CBL},
\label{eq:dtr2} \\
\mathcal{V}_{3}   :=   \mathcal{V}_{D(C)[1,\infty)}
= \mathcal{V}_{D(C)_c[1,\infty)} = \mathcal{W}_{B(C)_c[1,\infty)} = \mathcal{W}_{B(C)[1,\infty)} = \mathcal{W}_{B(C)} .
\label{eq:dtr3}
\end{gather}
The schematic \eqref{fig:simplify} is then simplied to Dixmier traces ($\mathcal{V}_1$), Connes--Dixmier traces ($\mathcal{V}_2$) and Ces\`{a}ro invariant Dixmier traces
($\mathcal{V}_3$):
\begin{equation} \label{fig:simplified_simplified}
\begin{array}{ccccccc}
\mathcal{V}_{\frac{1}{2}} & \supset &
\mathcal{V}_{1} & \supset &
\mathcal{V}_{2} & \supset &\mathcal{V}_{3} \\
& & \cup & & & &  \\
& & \mathcal{V}_{PDL[1,\infty)} & & & &
\end{array} .
\end{equation}
It is known $\mathcal{W}_{B(C)} \subsetneq \mathcal{W}_{CBL} \subsetneq \mathcal{W}_{BL}$, which is not trivial,
see~\cite{LSS} and \cite{SemSuk2009}.
Hence the inclusions are strict,
i.e.~$\mathcal{V}_1 \supsetneq \mathcal{V}_2 \supsetneq \mathcal{V}_3$ .
We do not know whether $\mathcal{V}_{\frac{1}{2}} \supset \mathcal{V}_{1}$ is strict, i.e.~anything is gained from the larger set of generalised limits. In some NCG papers Dixmier traces are def\/ined with vague reference to the type of generalised limit being used, relying
on the assumption of measurability (next section) to make the distinction irrelevant.  To our knowledge, the classes of functionals
in~\eqref{fig:simplified_simplified} represent the known relevant
distinctions in measuring the `log divergence of the trace' using the idea of Dixmier.

\subsection{Singular symmetric functionals} \label{sec:dtr:ssf}

It would be remiss to f\/inish a section on Dixmier traces without mentioning singular symmetric functionals.

Consider the following well-known formula for the canonical trace,
a special case of the Lidskii formula\footnote{There are Lidksii type formulations for the Dixmier trace, which for sake of brevity we will not discuss in this text, see~\cite{AS2005} and \cite{ssz2010}.  For spectral asymptotic formulae see \cite{Pin1995,FB}.}:
\begin{equation} \label{eq:lidskii_Tr}
\Tr(T) = \lim_{N \to \infty} \sum_{k=1}^N \mu_k(T)   , \qquad T > 0,
\end{equation}
which def\/ines a weight that is f\/inite on the ideal $\mathcal{L}^1$ of trace class operators.  The formula~\eqref{eq:lidskii_Tr} makes the
positivity and trace property evident (recall $\mu_k(U^*TU)=\mu_k(T)$
for all unitaries $U \in \mathcal{L}(H)$ and $T \in \mathcal{L}^\infty$), and makes it easy to see
$S\mapsto \Tr(ST)$, $S \in \mathcal{L}(H)$, $T \in \mathcal{L}^1$, is a positive continuous linear functional on $\mathcal{L}(H)$ by the same arguments as for the formula~\eqref{eq:dtr0}.

Linearity is not immediately evident from \eqref{eq:lidskii_Tr} as it was not for \eqref{eq:dtr0} (but is evident from the def\/inition of the trace
$\Tr(T) = \sum_m \langle h_m | T h_m \rangle$, $\{h_m\}_{m=1}^\infty$ an orthonormal basis of $H$, which will be of relevance in Section~\ref{sec:char}).

Both formula \eqref{eq:lidskii_Tr} and \eqref{eq:dtr0} have the general format
\begin{equation} \label{eq:symm_fun_0}
f( \{ \mu_n(T) \} ),
\end{equation}
where $f$ is a positive linear functional on the sequence space $\ell^1$
in the case of $\Tr$, and the sequence space
\[
m_{1,\infty} := \big\{ x= \{ x_n \}_{n=1}^\infty \in \ell^\infty \,| \, \sup_k \alpha(x^*)_k < \infty \big\} ,
\]
where $\alpha$ is the sequence from \eqref{eq:alpha_first}
and $x^*$ is the sequence $\{ |x_n| \}_{n=1}^\infty$ rearranged in nonincreasing order, in the case of $\Tr_\omega$.

A Calkin sequence space $j$ is a subset of $c_0$ such that
$x \in j$ and $y^* \leq x^*$ implies $y \in j$ (called the Calkin property) \cite[p.~18]{S}.  It is known that
a subset $\mathcal{J}$ of compact operators
is a two-sided ideal if and only if the sequence space
$j_{\mathcal{J}} = \{ x \in c_0 | x^* = \{ \mu_n(T) \}_{n=1}^\infty, T \in \mathcal{J} \}$
is a Calkin space~\cite{Calkin1941}.

A symmetric sequence space is a Calkin space $e$ with norm $\nm{\cdot}_e$ such that $y^* \leq x^*$ implies $\nm{y}_e \leq \nm{x}_e$ $\fa\, x \in e$.  A symmetric ideal of compact operators
$\mathcal{E}$ is a two sided ideal with a norm $\nm{\cdot}_\mathcal{E}$
such that $\nm{STV}_\mathcal{E} \leq \nm{S}\nm{T}_\mathcal{E} \nm{V}$, $\fa\, S,V \in \mathcal{L}(H)$, $T \in \mathcal{E}$.
It is known (actually, only recently when~\cite{KS2008_Creolle}
extended the f\/inite dimensional result of J.~von Neumann in~\cite{JvN1937}) that~$\mathcal{E}$ is a~symmetric ideal if
and only if $j_\mathcal{E}$ is a symmetric sequence space
and that $\nm{T}_\mathcal{E} = \nm{\{\mu_n(T)\}_{n=1}^\infty}_{j_{\epsilon}}$.  For example $\ell^1$ (resp.~$m_{1,\infty}$) is the symmetric
sequence space associated to the symmetric ideal~$\mathcal{L}^1$
(resp.~$\mathcal{M}_{1,\infty}$).

These correspondences allow us to form the equation \eqref{eq:symm_fun_0} for any symmetric ideal $\mathcal{E}$ and
using $0< f \in j_\mathcal{E}^*$.  The functional def\/ined on $\mathcal{E}_+$ by~\eqref{eq:symm_fun_0} is positive, unitarily invariant, continuous (in the symmetric norm $\nm{\cdot}_\mathcal{E}$), but not necessarily linear.  To obtain
linearity, and so derive a trace on $\mathcal{E}$, one introduces
the notions of symmetric and fully symmetric functionals.  If $e$ is a symmetric sequence space, then a functional $0 < f \in e^*$ is called symmetric
if $y^* \leq x^*$ implies $f(y^*) \leq f(x^*)$.
If $x, y \in e$, we write $y \prec\prec x$ (and say $y$ is submajorised by $x$) when
\[
\sum_{n=1}^N y_n^* \leq \sum_{n=1}^N x_n^*
\qquad \fa \, N \in \mathbb{N}.
\]
A symmetric sequence space $e$ is called fully symmetric if
$x \in e$ and $y \prec\prec x$ implies $y \in e$
(the Calkin property with rearrangment order replaced by Hardy--Littlewood--P{\'{o}}lya submajorisation~\cite{HLP1934}).
A functional $0 < f \in e^*$ on a fully symmetric sequence space
$e$ is called fully symmetric if $y \prec\prec x$ implies $f(y^*) \leq f(x^*)$.   Both $\ell^1$ and $m_{1,\infty}$ are fully symmetric sequence spaces.  The ideal $\mathcal{E}$ is called fully symmetric
if $j_{\mathcal{E}}$ is a fully symmetric sequence space.
Finally, a symmetric functional $0 < f \in e^*$ is called singular
if $f$ vanishes on all f\/inite sequences.

In \cite{DPSS, DPSSS1, DPSSS2}
and \cite[p.~77]{KS2008_Can}\footnote{These texts called symmetric (resp.~rearrangement invariant) the functionals which are here called fully symmetric (resp.~symmetric).}, it was shown (singular) symmetric linear
functionals
$0 < f \in j_{\mathcal{E}}^*$
correspond to (singular) positive linear unitarily
invariant functionals on fully symmetric ideals~$\mathcal{E}$ using the formula~\eqref{eq:symm_fun_0}.

Now, the relevant question is whether all Dixmier traces ($\mathcal{V}_1$) on $\mathcal{M}_{1,\infty}$ are constructed from singular symmetric functionals $0 < f \in m_{1,\infty}^*$.
The answer is yes.   In fact, we now know $f(\{ \mu_n(T) \})$ def\/ines a Dixmier trace if and only if $0 < f \in m_{1,\infty}^*$ is a singular \emph{fully} symmetric functional  (a result of the second author, N.~Kalton and A.~Sedaev, to appear). (Fully) symmetric singular functionals on (fully) symmetric ideals thus generalise the construction of Dixmier (other generalisations exist, see \cite{GI1995,AGPS1996}).  Another relevant question is whether
all positive linear unitarily invariant singular functionals on $\mathcal{M}_{1,\infty}$ constructed from singular symmetric functionals are Dixmier traces.  The answer is no.
There are singular symmetric functionals $0 < f \in m_{1,\infty}^*$ which are not fully symmetric,
yet $f(\{ \mu_n(T) \})$ def\/ines a singular positive linear unitarily invariant functional on $\mathcal{M}_{1,\infty}$ \cite[p.~79]{KS2008_Can}, i.e.~we now know
\begin{equation} \label{eq:traces_rho}
\mbox{there exist singular traces $0 < \rho \in \mathcal{M}_{1,\infty}^*$ that are not Dixmier traces.}
\end{equation}
In particular,  \cite[p.~79]{KS2008_Can} informs us that
there exists a non-trivial symmetric functional $0 < f_p \in m_{1,\infty}^*$ such that
$f_p(\{ \frac{1}{n} \}_{n=1}^\infty) = 0$
for the harmonic sequence $\{ \frac{1}{n} \}_{n=1}^\infty \in m_{1,\infty}$ (and hence it is not fully symmetric\footnote{All sequences in the unit ball of $m_{1,\infty}$ are submajorised by the harmonic sequence.  A non-trivial conti\-nuous functional on $m_{1,\infty}$ that vanishes on the harmonic sequence is not fully symmetric, else non-triviality is contradicted.}).
Setting $f = f_1 + f_p$,
where $f_1$ is any symmetric functional on $m_{1,\infty}$
for which $f_1(\{ \frac{1}{n} \}_{n=1}^\infty) = 1$,
yields a symmetric functional which cannot be fully symmetric.
The positive singular trace $\rho$ def\/ined by $\rho(T)= (f_1 + f_p)(\{\mu_n(T)\})$ is not a Dixmier trace, yet
$\rho(T) = 1$ for any positive compact operator $T$ with singular values given by the harmonic sequence (such as the inverse of the square root of the Laplacian on the 1-torus $\mathbb{R} / \mathbb{Z}$).

\section{Measurability} \label{sec:meas}

Measurability, as def\/ined by Connes~\cite[p.~308]{CN}, is the notion
that
\begin{equation} \label{eq:meas_0}
\Tr_\omega(T) = \Tr_{\omega'}(T) \qquad \fa \, \omega,\omega' \in B,
\end{equation}
where $B$ is some subset of generalised limits
which def\/ines a Dixmier trace (see the previous section).
Specif\/ically, Connes used the set $B = CDL_c[1,\infty)$ \cite[p.~305]{CN}.
A compact operator $T \in \mathcal{M}_{1,\infty}$ satisfying \eqref{eq:meas_0} is called $B$-measurable.

With reference to the noncommutative integral \eqref{eq:ncg_int},
one desires that the products $a \ds^{-n}$ are ($B$-)measurable,
so that the formula~\eqref{eq:ncg_int} is independent of the generalised limit used.  That the projections in $\A''$ (using the extension~\eqref{eq:ncg_ext}) have a unique value corresponding to the `inf\/initesimal' $\ds^{-n}$ parallels the idea that the Lebesgue measurable sets are those with a~unique `size' based on the metric.

Def\/ine the \emph{Tauberian} operators
\begin{equation} \label{eq:Taub}
\mathcal{T}^+ :=
\left\{ 0 < T \in \mathcal{M}_{1,\infty} \Big| \left\{ \frac{1}{\log(1+k)} \sum_{n=1}^k \mu_n(T) \right\}_{k=1}^\infty \in c \right\} .
\end{equation}
By construction, for a subset of generalised limits $B \subset S_\infty(\ell^\infty)$,
\[
\Tr_\omega(T) = \lim_{k \to \infty} \frac{1}{\log(1+k)} \sum_{n=1}^k \mu_n(T)
\qquad \fa \, \omega \in B , \quad T \in \mathcal{T}^+ .
\]
Hence the Tauberian operators are $B$-measurable for any $B \subset S_\infty(\ell^\infty)$.
We know from~\eqref{fig:simplified_simplified} in the previous section there are three (possibly f\/ive) main varieties of Dixmier trace. Then, \textit{a~priori}, there are three (possibly f\/ive)
versions of measurability. Def\/ine
\begin{equation} \label{eq:K}
\mathcal{K}_i :=
\{ T \in \mathcal{M}_{1,\infty} \, |\, \mathrm{f}(T)= \mathrm{const}
\ \fa \, \mathrm{f} \in \mathcal{V}_i \}   , \qquad
i=\frac{1}{2},1,2,3,PDL[1,\infty)
\end{equation}
and let $\mathcal{K}_i^+$ denote $0 < T \in \mathcal{K}_i$.
We concentrate on positive operators for the moment.
Then, from \eqref{fig:simplified_simplified}, we know that
\begin{equation} \label{eq:meas_seq_0}
\begin{array}{ccccccccc}
\mathcal{T}^+ & \subset & \mathcal{K}_{\frac{1}{2}}^+ & \subset &
\mathcal{K}_{1}^+ & \subset & \mathcal{K}_{2}^+ & \subset & \mathcal{K}_{3}^+ \\
& & & & \cap & & & &  \\
& & & & \mathcal{K}_{PDL[1,\infty)}^+ & & & &
\end{array}
\end{equation}
where the inclusions may or may not be strict.
In \cite{LSS}, our main result was that $\mathcal{K}_2^+ = \mathcal{T}^+$, identifying Connes' def\/inition of measurability in \cite{CN} for positive operators with the Tauberian operators.
Specif\/ically we introduced the notion of the Ces\`{a}ro limit property~\cite[p.~97]{LSS}
and any subset of Banach limits $B$ satisfying this property
has $B$-measurable positive operators (using the equivalence $\mathcal{V}_{DL[1,\infty)} = \mathcal{W}_{BL}$)
equivalent to Tauberian operators.  The set
$CBL$ of Ces\`aro--Banach limits satisfy the Ces\`aro limit property.  From the paper of E.~Semenov and the second author~\cite{SemSuk2009}, we now know that
$\mathcal{K}_{3}^+ \not=\mathcal{T}^+$.
From this, the schematic \eqref{eq:meas_seq_0} can be
reduced to
\begin{equation} \label{eq:meas_seq_1}
\begin{array}{ccccccccc}
\mathcal{T}^+ & = & \mathcal{K}_{\frac{1}{2}}^+ & = &
\mathcal{K}_{1}^+ & = & \mathcal{K}_{2}^+ & \subsetneq & \mathcal{K}_{3}^+ \\
& & & & \cap & & & &  \\
& & & & \mathcal{K}_{PDL[1,\infty)}^+ & & & &
\end{array} ,
\end{equation}
so that the set of positive measurable operators according to Connes--Dixmier traces (and any larger set of Dixmier traces)
is equivalent to the Tauberian operators, but the set of positive measurable operators for Ces\`{a}ro invariant Dixmier traces
is strictly larger than the Tauberian operators.
The Ces\`{a}ro limit property is suf\/f\/icient for equality with the Tauberian operators, it is an open question as to whether the property is a necessary condition (the set $B(C)$ of
Ces\`{a}ro invariant Banach limits does not satisfy the property
and $\mathcal{K}_{3}^+ \not=\mathcal{T}^+$).
Whether $\mathcal{K}_{PDL[1,\infty)}^+$ is equal to the Tauberian operators is also an open question.

It is interesting to note that some notions of measurability can exclude the Tauberian operators.  From the last section on symmetric functionals, since the set of Dixmier traces
corresponds to the set of fully symmetric functionals $0 < f \in m_{1,\infty}^*$,
the Tauberian operators $\mathcal{T}^+$ are the set of $0 < T \in \mathcal{M}_{1,\infty}$
such that the \eqref{eq:symm_fun_0} is independent of the fully symmetric functional~$f$.
Given~\eqref{eq:traces_rho}, we could try to widen this notion.
For example, let $\mathcal{V}_{\textrm{sym}}^+$ to be the set of all positive singular traces
$\rho$ on $\mathcal{M}_{1,\infty}$ such that $\rho(T) = f(\{ \mu_n(T) \}) \, \fa T > 0$
for some symmetric functional $0 < f \in m_{1,\infty}^*$.   Then set
\[
\mathcal{K}_{\textrm{sym}}^+ := \{ 0 < T \in \mathcal{M}_{1,\infty} \, | \,
\rho(T) = \mathrm{const} \ \fa \, \rho \in \mathcal{V}_{\textrm{sym}}^+ \} .
\]
We are informed by  \cite[p.~79]{KS2008_Can} that $\mathcal{T}^+ \not\subset \mathcal{K}_{\textrm{sym}}^+$.
So, widening the notion of measurable in this manner excludes some Tauberian operators.
We are indirectly informed by this result that there exist singular traces
on $\mathcal{M}_{1,\infty}$ that do not arise from
a generalised limiting procedure applied to the sequence
$\{ \frac{1}{\log(1+n)} \sum_{k=1}^n \mu_k(T) \}_{n=1}^\infty$.

In summary, the best way to avoid dif\/f\/iculties with the notion of measurability when one is using Dixmier traces to def\/ine the noncommutative integral \eqref{eq:ncg_int}
is to have the Tauberian condition:
\begin{equation} \label{cond:taub}
a^\frac{1}{2} \ds^{-n} a^\frac{1}{2} \in \mathcal{T}^+ , \qquad \fa \, 0 < a \in \A .
\end{equation}
However, we know of no general criteria involving the $*$-algebra
$\A$ and the selfadjoint opera\-tor~$D$ that demonstrates~\eqref{cond:taub}
besides direct examination of singular values.
Finding such a criteria is an open problem.
The linear span $\mathcal{T}$ of $\mathcal{T}^+$ is not an ideal
of compact operators, indeed we know $\mathcal{K}_i$, $i=1,2,3$,
is not an ideal for any class of Dixmier trace.

There are specif\/ic situations where \eqref{cond:taub} is implied without having to check singular values.  The functionals def\/ined by \eqref{eq:ncg_int} are all traces on the unital $*$-algebra $\A$.  If $\ds^{-n}$ is a Tauberian operator and $\Tr_\omega(\ds^{-n}) > 0$, then the traces from \eqref{eq:ncg_int} can be normalised by the same constant independent of $\omega$,
\begin{equation} \label{ncg_int_normalised}
\frac{\Tr_\omega(a\ds^{-n})}{\Tr_\omega(\ds^{-n})} .
\end{equation}
A unique trace state on $\A$ automatically implies the trace states in \eqref{ncg_int_normalised}
have the same value for all $\omega \in D_2$.  Hence, for all $a>0$, $a^{1/2} \ds^{-n} a^{1/2} \in \mathcal{T}^+$ and condition \eqref{cond:taub} is satisf\/ied.  The principle of a unique (singular) trace is also the construct for measurability of
pseudo-dif\/ferential operators of order $-n$.  So far as we know,
these situations of unique trace, besides very strong geometric conditions involving Hochschild cohomology and axioms for noncommutative manifolds \cite[p.~309]{CN}, \cite[p.~160]{C4}
provide the only handle on~\eqref{cond:taub}.

The unsatisfactory feature of this situation, from a measure-theoretic point of view, is that a~unique
trace often enables more general traces than Dixmier traces on compact operators to be used in \eqref{eq:ncg_int} to obtain the same result (see the sections below).
Where is the def\/initive rationale that Dixmier traces
should solely appear in connection with the noncommutative integral?

\subsection{Decomposition of measurable operators}

We concentrated in the previous section on identifying those
\emph{positive} compact operators for which $\Tr_\omega(T) = \mathrm{const}$ for various $\Tr_\omega$ in the classes $\mathcal{V}_1$, $\mathcal{V}_2$, $\mathcal{V}_3$.
We know that general $T\in\mathcal{M}_{1,\infty}$ polarise, i.e.
\[
T = T_1 - T_2 + iT_3 - iT_4  , \qquad 0 < T_j \in \mathcal{M}_{1,\infty}  , \qquad j=1,2,3,4 ,
\]
and that
\begin{equation} \label{eq:dtr_polar}
\Tr_\omega(T) =
\omega \left( \left\{ \frac{1}{\log(1+N)} \sum_{k=1}^N
\mu_k(T_1) - \mu_k(T_2) +  i\mu_k(T_3) - i\mu_k(T_4) \right\}_{k=1}^\infty \right).
\end{equation}
For convenience, denote
\[
\tilde{\mu}_k(T) := \mu_k(T_1) - \mu_k(T_2) +  i\mu_k(T_3) - i\mu_k(T_4)  , \qquad T \in \mathcal{L}^\infty(H) .
\]
It is evident if $T \in \mathcal{T} := \mathrm{span}_{\mathbb{C}}(\mathcal{T}^+)$ (i.e.~$T_1,T_2,T_3,T_4 \in \mathcal{T}^+$), or
\[
T \in \tilde{\mathcal{T}} :=
\left\{ T \in \mathcal{M}_{1,\infty} \, \Big|\, \left\{ \frac{1}{\log(1+N)} \sum_{k=1}^N \tilde{\mu}_k(T) \right\}_{k=1}^\infty \in c \right\} \supset \mathcal{T} ,
\]
then
\[
\Tr_\omega(T) = \lim_{N \to \infty}  \frac{1}{\log(1+N)} \sum_{k=1}^N
\tilde{\mu}_k(T) .
\]
Hence $\Tr_\omega(T) = \mathrm{const}$ for $T \in \tilde{\mathcal{T}}$.  Does the reverse
implication hold for the various classes $\Tr_\omega \in \mathcal{V}_1,\mathcal{V}_2,\mathcal{V}_3$?

Recent work by the second author, N.~Kalton, and A.~Sedaev, to appear, indicates that the proposition
\begin{equation} \label{eq:meas_U_0}
\mathrm{f}(T) = \mathrm{const} \  \fa \,\mathrm{f} \in \mathcal{V}_1 \ \Rightarrow \
\left\{ \frac{1}{\log(1+N)} \sum_{k=1}^N \tilde{\mu}_k(T) \right\}_{k=1}^\infty \in c
\end{equation}
is~\emph{false}.  In the next section on symmetric subideals
we brief\/ly describe the strictness of the inclusions in the following schematic:
\begin{equation} \label{eq:meas_seq_2}
\begin{array}{ccccccc}
\tilde{\mathcal{T}} & \subsetneq &
\mathcal{K}_{1} & \subsetneq & \mathcal{K}_{2} & \subsetneq & \mathcal{K}_{3}
\end{array} .
\end{equation}
The place of $\mathcal{K}_{\frac{1}{2}}$ (and
$\mathcal{K}_{PDL[1,\infty)}$) is open.

The strictness of the inclusion $\tilde{\mathcal{T}} \subsetneq \mathcal{K}_{1}$ and the result \eqref{eq:meas_seq_1}
for positive elements implies
\[
\mathrm{span}_{\mathbb{C}}(\mathcal{K}_i^+) \subsetneq \mathcal{K}_{i}   , \qquad i=1,2 .
\]
It is evident from $\mathrm{f}(T^*) = \overline{\mathrm{f}(T)}$, $\mathrm{f} \in \mathcal{V}_i$,
$i = \frac{1}{2},1,2, PDL[1,\infty)$,
that the sets $\mathcal{K}_i$ are closed under the $*$-operation.
It follows that $T \in \mathcal{K}_i \not\Rightarrow |T| \in \mathcal{K}_i$, $i=1,2$.

As a consequence, we lack a basic decomposition result
of integration theory. Suppose $T=T^*$ with the standard unique
decomposition into positive and negative operators
\begin{equation} \label{eq:pos_decomp}
T  = T_+ - T_-   , \qquad T_{\pm} := \frac{1}{2}(|T| \pm T) .
\end{equation}
Then,
\[
T=T^* \in \mathcal{K}_i  \ \not\Rightarrow \
T_+,T_- \in \mathcal{K}_i   , \qquad i=1,2 .
\]
Compare this with the result for a Borel measure space $(X,\mu)$
(also true for the noncommutative variant of $\tau$-measurable operators on $(\mathcal{N},\tau)$ for a semif\/inite von Nuemann algebra $\mathcal{N}$ with faithful normal semif\/inite trace $\tau$)
\[
f=\overline{f} \in L^0(X,\mu) \Rightarrow
\exists! \, 0 < f_+,f_- \in L^0(X,\mu) \qquad \mbox{s.t.} \ f = f_+-f_-,
\]
where $L^0(X,\mu)$ denotes the set of (equivalence classes) of $\mu$-measurable functions on $X$.

\subsection{Closed symmetric subideals} \label{sec:dtr:subideals}

Def\/ine symmetric subspaces of the symmetric sequence space $m_{1,\infty}$ by
\begin{gather}
u_{1,\infty}   =   \left\{ x \in m_{1,\infty} \, | \, x^*(n) = o\left(\frac{\log n}{n}\right) \right\}, \label{def:u} \\
\ell^{1,w}   =   \left\{ x \in m_{1,\infty} \, | \, x^*(n) = O\left(\frac{1}{n}\right) \right\} . \label{def:weak_l1}
\end{gather}
The ideal of compact operators $\mathcal{U}_{1,\infty}$
and $\mathcal{L}^{1,w}$, def\/ined by
\[ u_{1,\infty} =j_{\mathcal{U}_{1,\infty}}   , \qquad
\ell^{1,w} =j_{\mathcal{L}^{1,w}} ,
\]
are quasi-Banach symmetric ideals.
Notice also
\[
\mathcal{L}^{1,w} \subsetneq \mathcal{U}_{1,\infty}
\subsetneq \mathcal{M}_{1,\infty} .
\]
Def\/ine, as before,
\[
\tilde{\mathcal{T}} :=
\left\{ T \in \mathcal{M}_{1,\infty} \Big| \left\{ \frac{1}{\log(1+N)} \sum_{k=1}^N \tilde{\mu}_k(T) \right\}_{k=1}^\infty \in c \right\} .
\]
New results of N.~Kalton, A.~Sedaev, and the second author,
state that
\begin{equation} \label{eq:u_1}
\tilde{\mathcal{T}} \cap \mathcal{U}_{1,\infty} =
\mathcal{K}_1 \cap \mathcal{U}_{1,\infty}
\end{equation}
and that $\mathcal{U}_{1,\infty}$ is the maximal symmetric subideal
$\mathcal{E} \subset \mathcal{M}_{1,\infty}$
for which \eqref{eq:u_1} holds.
Thus,
\[
\tilde{\mathcal{T}} = \tilde{\mathcal{T}} \cap \mathcal{M}_{1,\infty}
\not= \mathcal{K}_1 \cap \mathcal{M}_{1,\infty} = \mathcal{K}_1 .
\]
Their results also show
\begin{equation} \label{eq:u_2a}
\tilde{\mathcal{T}} \cap \mathcal{U}_{1,\infty} \not=
\mathcal{K}_2 \cap \mathcal{U}_{1,\infty} .
\end{equation}
Hence $\mathcal{K}_1 \not= \mathcal{K}_2$.  It is an open question as to the existence of a maximal closed symmetric subideal $\mathcal{U}_2$ such that $\tilde{\mathcal{T}} \cap \mathcal{U}_{2} = \mathcal{K}_2 \cap \mathcal{U}_{2}$.  Is $\mathcal{L}^{1,w} \subseteq \mathcal{U}_2$?

\section{Origin of the noncommutative integral} \label{sec:origin}

The origin of the formula \eqref{eq:ncg_int} as the integral
in noncommutative geometry lies in
Connes' trace theorem \cite[Theorem~1]{C3}, \cite[p.~293]{GBVF}:

\begin{theorem}[Connes' trace theorem] \label{thm:2}
Let $M$ be a compact $n$-dimensional manifold, $\mathcal{E}$~a~complex vector bundle
on $M$, and $P$ a pseudo-differential operator of order $-n$ acting on sections of $\mathcal{E}$.
Then the corresponding operator $P$ in $H=L^2(M,\mathcal{E})$ belongs to $\mathcal{M}_{1,\infty}(H)$ and one has:
\begin{equation} \label{eq:cntr}
\Tr_\omega(P) =  \frac{1}{n} \mathrm{Res}(P)
\end{equation}
for any $\Tr_\omega \in \mathcal{V}_1$.
\end{theorem}

Here $\mathrm{Res}$ is the restriction of the
Adler--Manin--Wodzicki residue to pseudo-dif\/ferential operators of order $-n$ \cite{Wod, C3}.  The monograph \cite[\S~7]{GBVF} provides
a detailed introduction to the residue $\mathrm{Res}$
and the origin and framework of the formula.  We recall that,
for classical pseudo-dif\/ferential operators of order $-n$ acting on the trivial line bundle ($H=L^2(M)$),
\[
\mathrm{Res}(P) = \frac{1}{(2\pi)^n} \int_{S^1M} \sigma_{-n}(P)(x,\xi) dx d\xi,
\]
where $dx d\xi$ is the Liouville measure (volume form) on the cosphere bundle $S^1M$
(with f\/ibres $S^1M_x \cong \{ x \} \times S^{n-1}$ over $x \in M$ where $S^{n}$ is the
unit $n$-sphere)
and $\sigma_{-n}$ is the principal symbol of~$P$.
On the left hand side of \eqref{eq:cntr} one should read
for $P$ the unique extension of the negative order pseudo-dif\/ferential operator
$C^\infty(M) \to C^\infty(M)$ to a compact operator $L^2(M) \to L^2(M)$~\cite{Shubin2003}.

Let $\Delta$ denote the Hodge Laplacian on $M$
and $P_a \in B(L^2(M))$ denote the extension of the zeroth order pseudo-dif\/ferential operator associated to a function $a \in C^\infty(S^1M)$.  Since
\[
\sigma_{-n}(P_a(1+\Delta)^{-n/2})|_{S^1M}
= a
\]
formula \eqref{eq:cntr} provides
\begin{equation} \label{eq:cntr2}
\Tr_\omega\big(P_a(1+\Delta)^{-n/2}\big) =  \frac{1}{n(2\pi)^n} \int_{S^1M} a(x,\xi) dx d\xi .
\end{equation}
From the inclusion $a \in C^\infty(M) \hookrightarrow C^\infty(S^1M)$ the operator $P_a$, $a \in C^\infty(M)$, is the multiplier by the function $a \in C^\infty(M)$ and \cite[Corollary~7.21]{GBVF}, \cite[\S~1.1]{FB} or \cite[p.~98]{Landi},
\begin{equation} \label{eq:integral}
\Tr_\omega(P_a T_\Delta) = \frac{\mathrm{Vol}(S^{n-1})}{n(2\pi)^n} \int_{M} a(x) dx  , \qquad a \in C^\infty(M),
\end{equation}
where we set
\[
T_\Delta := (1+\Delta)^{-n/2} \in \mathcal{M}_{1,\infty}.
\]
Equation \eqref{eq:integral} originated the use of the term integral for the expression \eqref{eq:ncg_int}.

Why is a singular trace needed in formula \eqref{eq:cntr}?
By the asymptotic expansion of classical symbols $P - P_0 \equiv \sigma_{-n}(P)$
for a classical pseudo-dif\/ferential operator $P$ of order $-n$
where $P_0$ is of order $< -n$ and $\sigma_{-n}(P)$ is the principal symbol of $P$ \cite[XVIII]{HorIII}.
All pseudo-dif\/ferential operators of order $< -n$ on a
compact $n$-dimensional manifold have trace class extensions~\cite{Shubin2003}.
Hence any singular functional $\mathrm{f} \in \mathcal{M}_{1,\infty}$ applied to $P$, i.e.~$P \mapsto \mathrm{f}(P)$, since it vanishes on $\mathcal{L}^1$, is equivalent to a functional on symbols of order $-n$.
Moreover the vanishing on $\mathcal{L}^1$ permits that the statement
need only be proved locally.  Hence $M$ may be taken to be homogeneous.
If $\mathrm{f}$ is a positive trace this functional on symbols of order~$-n$
is equivalent to an invariant measure on~$S^1M$.  Hence, up to a constant, it
identif\/ies with the Liouville integral.

The above argument, from the proof of
Connes' trace theorem in~\cite{C3}, implies that any
positive singular trace $0 < \rho \in \mathcal{M}_{1,\infty}^*$
can be substituted for the Dixmier
trace in \eqref{eq:cntr}.  Therefore we may write, using the above notation:

\begin{theorem}[Connes' trace theorem] \label{thm:3}
Let $M$ be a compact $n$-dimensional manifold, $\mathcal{E}$~a~complex vector bundle
on $M$, and $P$ a pseudo-differential operator of order~$-n$ acting on sections of $\mathcal{E}$.
Then the corresponding operator $P$ in $H=L^2(M,\mathcal{E})$ belongs to $\mathcal{M}_{1,\infty}(H)$ and one has:
\begin{equation} \label{eq:cntr3}
\rho(P) =  \frac{\rho(T_\Delta)(2\pi)^n}{\mathrm{Vol}(S^1M)} \mathrm{Res}(P)
\end{equation}
for any singular trace $0 < \rho \in \mathcal{M}_{1,\infty}(H)^*$ .
\end{theorem}

Singular traces $0 < \rho \in \mathcal{M}_{1,\infty}(L^2(M,\mathcal{E}))^*$, which are not Dixmier traces
and such that $\rho(T_\Delta) > 0$, exist by a variant of
the argument after \eqref{eq:traces_rho}.  Correspondingly
\begin{equation} \label{eq:integral2}
\rho(P_a T_\Delta) = \frac{\rho(T_\Delta)}{\mathrm{Vol}(M)} \int_{M} a(x) dx  , \qquad a \in C^\infty(M)
\end{equation}
provides the integral on a manifold for smooth functions
up to a non-zero constant.

The result \eqref{eq:integral2} raises the question about whether
the class of Dixmier traces ($\mathcal{V}_1$) is too restrictive
in the def\/inition of the noncommutative integral
\eqref{eq:ncg_int}.  If we def\/ine
\begin{equation} \label{def:V_sing}
\mathcal{V}_{\mathrm{sing}}^+ := \big\{ 0 < \rho \in \mathcal{M}_{1,\infty}^* \,|\,
\rho \mbox{ \ is a singular trace} \big\}
\end{equation}
and then def\/ine as the `integral'
\begin{equation} \label{eq:ncg_int2}
\rho(a\ds^{-n})  , \qquad a \in \A , \qquad \rho \in \mathcal{V}_{\mathrm{sing}}^+ , \qquad \rho(\ds^{-n}) > 0
\end{equation}
we have the same claim to formulas \eqref{eq:cntr} and \eqref{eq:integral}.  Of course, there are a great many formulas involving Dixmier traces besides \eqref{eq:cntr} and \eqref{eq:integral}, particularly in the development of (local) index theory in noncommutative geometry \cite{CM,FB,AF, CPRS_Hoch2004, CPRSI2006,CPRSII2006,CaPoSu,Surveys2000Hig}.
There are geometric consequences of the formulation involving the Dixmier trace and associations to pseudo-dif\/ferential operators \cite{C3, Ponge,GI2000, Hankel2009}, including extensions of Connes' trace theorem \cite{Nicola2008, Nest97dix, Ponge}.  There are also many uses of the Dixmier trace in genuine noncommutative examples, for example
in the geometry of totally disconnected sets \cite{Bell2009,Lap1997, GI2003}.

Are there any of these results for which the submajorisation
property of the Dixmier trace is a~necessity (recall Dixmier traces correspond to the fully symmetric singular functionals on~$m_{1,\infty}$
from Section~\ref{sec:dtr:ssf})?  Several current proofs
depend on submajorisation, e.g.~residue and heat kernel
results are dependent on submajorisation (cf.~proofs in \cite{CRSS} and \cite{ssz2010}).
It is open whether this dependence is a necessity for the results or an artifact of the proofs.

\subsection{Normal extension}

To make things simpler, we dispense with the constants for the moment.  Def\/ine the noncommutative integral as
\begin{equation} \label{eq:ncg_state}
\Phi_\rho(a) = \rho(a\ds^{-n})  , \qquad a \in \A , \qquad \rho \in \mathcal{V}_{\mathrm{sing}}^+  , \qquad  \rho(\ds^{-n}) = 1 .
\end{equation}
By the discussion in the introduction,
the same formula $\Phi_\rho$ def\/ines a state of $A$ and $\A''$.

Let $M$ be a compact $n$-dimensional Riemannian manifold, $\A = C^\infty(M)$ acting by multipliers on $L^2(M)$, and $d+d^*$
the Hodge--Dirac operator acting (densely) on square integrable sections of the exterior bundle.  Then
\[
\langle d+d^* \rangle^{-n} = (1+\Delta)^{-n/2} =: T_\Delta
\]
which we take as acting as a compact operator on $L^2(M)$.
By \eqref{eq:integral2}
\begin{equation} \label{eq:integral3}
\Phi_\rho(f) = \int_M f d \hat{x}  , \qquad f \in C^\infty(M)
\end{equation}
where $d\hat{x}$ is the normalised volume form.
Recall $\Phi_\rho$ is a state on $L^\infty(M)$.  Hence, by uniform continuity
of both the left hand side and right hand side of~\eqref{eq:integral3}
\[
\Phi_\rho(f) = \int_M f d \hat{x}  , \qquad f \in C(M) ,
\]
which identif\/ies $\Phi_\rho|_{C(M)}$ with the normalised
integral on $C(M)$.   Hence $\Phi_\rho|_{C(M)} \in C(M)^*$ has an extension to $L^\infty(M)$ given by the normalised integral
\[
\phi(f) = \int_M f d \hat{x}  , \qquad f \in L^\infty(M) .
\]
If one considers for the moment a homogeneous compact manifold $M$,
by inspection $\Phi_\rho \in L^\infty(M)^*$ is a invariant state,
and one might think it is obvious that we must have
$\Phi_\rho = \phi$ since the volume form provides the unique invariant measure class on $M$.  While this last phrase is true, there is no reason \textit{a~priori} why the state
$\Phi_\rho \in L^\infty(M)^*$ corresponds to a measure.
We know that there can be an inf\/initude of invariant states
on $L^\infty(M)$ which all agree with the Lebesgue integral on $C(M)$~\cite[Theorem~3.4]{Rudin} (and f\/irst Baire class functions\footnote{We are indebted to B.~de Pagter for pointing this out and bringing Rudin's paper to our attention.  We also thank P.~Dodds for additional explanation.}).

It was not until recently that an elementary argument\footnote{The generalisation of this argument is discussed in our f\/inal section, Section~\ref{sec:char}.} was found
proving $\Phi_\rho = \phi$ for the $n$-torus.
An earlier claim \cite[p.~297]{GBVF},
applied monotone convergence to both sides of \eqref{eq:integral3} to extend from $C^\infty$-functions to $L^\infty$-functions.
This method of proof is circular, monotone convergence can be applied to the right hand side, since the integral is a normal linear function on~$L^\infty(M)$.
To apply monotone convergence to the left hand side it must be known $\Phi_\rho \in L^\infty(M)_*$.

In \cite{LDS} we used a dif\/ferent method than that for torii
to show the following result.  We remark that~\cite{LDS} concentrated only on Dixmier traces.  We state the result
for general singular traces here.

\begin{theorem} \label{thm:result2}
Let $M$ be a compact $n$-dimensional Riemannian manifold, $\Delta$
the Hodge--Lap\-la\-cian on $L^2(M)$, and $T_\Delta$ as above.
Then $f T_\Delta \in \mathcal{M}_{1,\infty}(L^2(M))$ if and only if $f \in L^2(M)$ and
\begin{equation} \label{eq:integral4}
\Phi_{\rho}(f) := \rho(f T_\Delta) = \int_{M} f d\hat{x} , \qquad \fa \, f \in L^2(M)
\end{equation}
for any singular trace $0 < \rho \in \mathcal{M}_{1,\infty}(L^2(M))^*$ such that $\rho(T_\Delta) = 1$.
\end{theorem}
\begin{proof}
We adapt only those parts of  \cite[Theorem 2.5]{LDS} concerning~\eqref{eq:integral4}.  Def\/ine the sequence space
\begin{equation} \label{eq:m0}
m_0  =  \overline{\ell^1} \ \ \mbox{in the norm of} \ \ m_{1,\infty} \end{equation}
and the ideal $\mathcal{M}_0$ given by $m_0 =j_{\mathcal{M}_0}$.
The Riesz seminorm on $\mathcal{M}_{1,\infty}$ is def\/ined by
\begin{equation} \label{eq:norm0}
\nm{T}_0 := \inf_{Q \in \mathcal{M}_0} \nm{T-Q}_{1,\infty} .
\end{equation}
By construction every singular functional
$0 < \rho \in \mathcal{M}_{1,\infty}^*$ vanishes on $\mathcal{M}_0$.
Hence
\[
|\rho(T)| \leq \nm{\rho} \nm{T}_0  , \qquad T \in \mathcal{M}_{1,\infty} .
\]
The proof in \cite{LDS} involved, f\/irstly, showing that
\[
\nm{f T_\Delta}_0 \leq C \nm{f}_2 \nm{T_\Delta}_{1,\infty}
\]
for a constant $C>0$ (achieved by  \cite[Corollary~4.5]{LDS}  and \cite[Examples~4.6,~4.7]{LDS}).
Hence, from the preceeding two equations, $f T_\Delta \in \mathcal{M}_{1,\infty}$ and
\[
|\rho(fT_\Delta)| \leq C \nm{\rho} \nm{f}_2 \nm{T_\Delta}_{1,\infty}
\]
for all $f \in C^\infty(M)$.  Moreover, if $C^\infty(M) \ni f_n \to f \in L^2(M)$ in the $L^2$-norm then $fT_\Delta$ is bounded,
indeed $f T_\Delta \in \mathcal{M}_{1,\infty}$,
and
\[
|\rho((f-f_n)T_\Delta)| \leq C \nm{\rho} \nm{f-f_n}_2 \nm{T_\Delta}_{1,\infty} \to 0.
\]
The result~\eqref{eq:integral4} now follows from Connes' trace theorem (specif\/ically,
from~\eqref{eq:integral3}).
\end{proof}

As a corollary,
\[
\Phi_\rho(f) = \int_M f d \hat{x}  , \qquad f \in L^\infty(M)
\]
and we have an example where
\[
\Phi_\rho \in L^\infty(M)_*
\]
(i.e.~the noncommutative integral is normal).

\subsection{Tools for the noncommutative integral}

As mentioned in the introduction, when
\begin{equation} \label{eq:ncg_state_normal}
\Phi_\rho(a) = \rho(a\ds^{-n})  , \qquad a \in \A'', \qquad \rho \in \mathcal{V}_{\mathrm{sing}}^+
  , \qquad \rho(\ds^{-n}) = 1
\end{equation}
def\/ines a normal state on the von Neumann algebra $\A'' \subset \mathcal{L}(H)$, there are
Radon--Nikodym theorems, dominated convergence theorems, etc.~for this state.

Very little is known about equivalent measure-theoretic results
when \eqref{eq:ncg_state_normal} is not known to be normal
(and hence the emphasis at present centres on characterising normality in \eqref{eq:ncg_state_normal}, see Section~\ref{sec:char}).

Def\/ine the classical weak $\ell^p$-spaces (or equivalently the
fully symmetric Lorentz sequence spaces $\ell^{p,\infty}$):
\begin{gather}
\ell^{p,w}    =   \left\{ x \in c_0 \,|\, x^*(n) = O\left(n^{-\frac{1}{p}}\right) \right\} , \qquad 1 \leq p < \infty. \label{def:weak_lp}
\end{gather}
Denote the corresponding closed fully symmetric ideals of compact operators by $\mathcal{L}^{p,w}$ (i.e.~$\ell^{p,w} = j_{\mathcal{L}^{p,w}}$).  As is classically known,
$\ell^{p,w}$ is the $p$-convexif\/ication of $\ell^{1,w}$, $1 \leq p < \infty$~\cite{LTza1979}.  Similarly, from~\cite{CRSS} or~\cite{KS2008_Creolle}, $\mathcal{L}^{p,w}$ is the $p$-convexif\/ication of $\mathcal{L}^{1,w}$.  The $p$-convexif\/ication $\mathcal{M}_{p,\infty}$ of the ideal~$\mathcal{M}_{1,\infty}$, denoted $\mathcal{Z}_p$ in \cite{CRSS}, is strictly larger that
$\mathcal{L}^{p,w}$
(see~\cite{CRSS})\footnote{Note that the $(p,\infty)$-summable condition for K-cycles (spectral triples) is sometimes stated as $\ds^{-1} \in \mathcal{L}^{p,w}$, $1\leq p < \infty$ \cite[p.~546]{CN}, \cite[p.~159]{C4}
(equivalent to $\ds^{-p} \in \mathcal{L}^{1,w}$)
and sometimes as $\ds^{-1} \in \mathcal{M}_{1,\infty}$ ($p=1$)
and $\ds^{-1} \in \mathcal{L}^{p,w}$ ($p>1$) \cite[IV.2]{CN} (not equivalent to $\ds^{-p} \in \mathcal{L}^{1,w}$
or $\ds^{-p} \in \mathcal{M}_{1,\infty}$).  To distinguish the strictly smaller ideal $\mathcal{L}^{1,w}$
from the usually quoted domain of the Dixmier trace (the ideal $\mathcal{M}_{1,\infty}$) we have avoided the notation
$\mathcal{L}^{1,\infty}$ for $\mathcal{M}_{1,\infty}$.}.

{\sloppy We do know, from F.~Cipriani, D.~Guido, and S.~Scarlatti~\cite{Cip1996}, and later Guido and T.~Isola~\cite{GI2003}, that Dixmier traces satisfy a H\"{older} inequality.
The result of  \cite[Theorem~5.1]{GI2003} applies to
all singular traces $\rho$ used in~\eqref{eq:ncg_state_normal}.

}

\begin{theorem}[\protect{\cite[Lemma~1.4]{Cip1996}, \cite[Theorem~5.1]{GI2003},
\cite[Lemma~6.2(i)]{CRSS}}]
For $T \in \mathcal{M}_{p,\infty}$ and $V \in \mathcal{M}_{q,\infty}$, $p,q>1$, $p^{-1}+q^{-1}=1$, and a singular trace
$0 < \rho \in \mathcal{M}_{1,\infty}^*$, we have
$TV \in \mathcal{M}_{1,\infty}$ and the
H\"{o}lder inequality:
\[
|\rho(TV)| \leq \rho(|TV|) \leq C_p \rho(|T|^p)^\frac{1}{p} \rho(|V|^q)^\frac{1}{q},
\]
where $C_p :=1$ if $\rho \in \mathcal{V}_1$
and $C_p := 1 + 2 \frac{\sqrt{p-1}}{p}$ otherwise.
The inequality also holds for $T \in \mathcal{M}_{1,\infty}$
and $V \in \mathcal{L}(H)$:
\[
|\rho(TV)| \leq \rho(|TV|) \leq \nm{V} \rho(|T|)
\]
$($the $p=1$, $q=\infty$ case$)$.
\end{theorem}

The H\"{o}lder inequality is used to prove that
$\rho(a\ds^{-n})$ is a~trace on
$\A$ when $(\A,H,D)$ is a~spectral triple and $\rho \in \mathcal{V}_1$ \cite{Cip1996},
\cite[Theorem~6.1]{CRSS}.  In \cite{GI2003}
it was noted the same proof applies to any singular trace.

\begin{corollary} \label{cor:trace}
Let $(\A,H,D)$ be a~spectral triple such that $\ds^{-n} \in \mathcal{M}_{1,\infty}$.  Then \eqref{eq:ncg_int2}
$($resp.~{\rm \eqref{eq:ncg_state_normal})}
defines a finite positive trace $($resp.~trace state$)$ on~$\A$.
\end{corollary}

\section{Zeta functions and heat kernels} \label{sec:res}

Zeta functions and heat kernel asymptotics are
alternative ways to measure `the log divergence of the trace'.

The Wodzicki residue $\mathrm{Res}$ on classical pseudo-dif\/ferential
operators of order $-n$ on a $n$-dimensional compact Riemannian manifold $M$ derives its name as the noncommutative residue
from the zeta function formulation of the residue for positive elliptic
pseudo-dif\/ferential operators~\cite{Wod, Wod1987}.

Explicitely, if $0 < Q \in \mathrm{Op}^{d}_{\mathrm{cl}}(M)$,
$d > 0$,
is elliptic, then $Q^{-s}$, $s>n/d$, is trace class \cite{Seeley1967,Shubin2003,HorIII}.
It is known the function
\[
\zeta_Q(s) := \Tr(Q^{-s}) , \qquad s>n/d  ,
\]
called the zeta function and initially introduced for the
Laplacian~\cite{MinPle1949} has a meromorphic continuation
with simple pole at $s=n/d$ \cite{Seeley1967,Hig2004}
and from \cite{Wod}
\begin{equation} \label{eq:res_classic}
\mathrm{res}_{s=n/d} \zeta_Q(s)
= \lim_{s \to (n/d)^+} (s-n/d)\Tr(Q^{-s})
= -\frac{1}{d} \mathrm{Res}(Q^{-n/d}).
\end{equation}
The Wodzicki residue has a similar identif\/ication involving
the heat kernel operator $e^{-tQ}$, $t > 0$.
The heat kernel operator is so named since the kernel $K(t,x,y) \in C^\infty((0,\infty) \times M \times M)$
associated to the trace-class family $e^{-t\Delta}$, $t > 0$,
$\Delta$ the Hodge Laplacian,
\[
(e^{-t\Delta}h)(x) = \int_{M}  K(t,x,y) h(y) dy  , \qquad h \in L^2(M)
\]
is a solution to the (local) heat equation \cite{MinPle1949}
\[
(\partial_t - \Delta_x)K(t, x, y) = 0  , \qquad \lim_{t \to 0^+} K(t, x, y) = \delta_y(x),
\]
where $\delta$ is the Dirac delta function and the limit is in the weak sense.  The zeta function and the trace of the heat kernel operator are related by the Mellin transform
\[
\Gamma(s) \zeta_Q(s) = \int_0^\infty t^{s-1} \Tr(e^{-tQ}) dt .
\]
Using the asymptotic expansion of the heat kernel operator \cite{MinPle1949, Gilkey1995, BGV1992}
the Wodzicki residue is associated to the heat kernel by
\begin{equation} \label{eq:htk_classic}
\lim_{t \to 0^+} t^{n/d} \Tr(e^{-tQ})
= \Gamma\left(\frac{n}{d}\right) \mathrm{res}_{s=n/d} \zeta_Q(s)
= -\frac{1}{d}\Gamma\left(\frac{n}{d}\right) \mathrm{Res}(Q^{-n/d}) .
\end{equation}
From Connes' trace theorem (Theorem~\ref{thm:2} and~\ref{thm:3})
the Dixmier trace of $\ds^{-n}$ has claim to be the noncommutative version of the
Wodzicki residue $\mathrm{Res}(Q^{-n})$, $0 < Q \in \mathrm{Op}^{1}_{\mathrm{cl}}(M)$.
Do residue and heat kernel formulas similar to~\eqref{eq:res_classic}
and~\eqref{eq:htk_classic} hold for Dixmier traces?
Sections~\ref{sec:res.1} and~\ref{sec:res.2} list what is known.

We detail only the latest known identif\/ications
between Dixmier traces and residues and heat kernels.
For applications of the residue and heat kernel to index formulations
in NCG see \cite{CM, AF, CPRSI2006,CPRSII2006, Surveys2000Hig}.

\subsection{Residues of zeta functions} \label{sec:res.1}

A.~Connes introduced the association between a generalised zeta function,
\begin{equation} \label{eq:zeta}
\zeta_T(s) := \Tr(T^s) = \sum_{n=1}^\infty \mu_{n}(T)^s , \qquad 0 < T \in \mathcal{M}_{1,\infty}
\end{equation}
and the calculation of a Dixmier trace with the result that
\begin{equation} \label{eq:res_orig}
\lim_{s \to 1^+} (s-1) \zeta_T(s) = \lim_{N \to \infty} \frac{1}{\log(1+N)} \sum_{n=1}^{N} \mu_{n}(T)
\end{equation}
if either limit exists \cite[p.~306]{CN}.
Meromorphicity of the formula~\eqref{eq:zeta} is discussed in~\cite{CM, Connes_GSPV,Hig2004}.
Generalisations of \eqref{eq:res_orig} appeared in~\cite{CPS} and later~\cite{CRSS}.

In \cite{LS2} we translated the results  \cite[Theorem 4.11]{CRSS} and  \cite[Theorem 3.8]{CPS} to $\ell^\infty$, see Theorem~\ref{thm:resPA} and Corollary~\ref{cor:resA} below.  Recall
$\mathcal{T}^+ = \mathcal{K}_1^+ = \mathcal{K}_2^+$ \eqref{eq:meas_seq_1}, so
\begin{equation} \label{eq:res0}
\Tr_\omega(T) = \lim_{s \to 1^+} (s-1) \zeta_T(s)
\end{equation}
calculates the Dixmier or Connes--Dixmier trace of any (Dixmier or Connes--Dixmier) measurable positive operator $0 < T \in \mathcal{M}_{1,\infty}$ as the residue at $s=1$ of the zeta function $\zeta_T$.  What about general $0 < T \in \mathcal{M}_{1,\infty}$?

The zeta function of a positive compact operator $T$ given by~\eqref{eq:zeta} relies on the assumption that there exists some $s_0$ for
which $T^s$ is trace class $s > s_0$ (equally $T \in \mathcal{L}^s$
for $s > s_0$).
For convenience we assume $s_0=1$.  The space
of compact operators for which the zeta function exists
and $(s-1)\zeta_T(s) \in L^\infty((1,2])$
was studied in~\cite{CRSS}.  Def\/ine the norm
\[
\nm{T}_{\mathcal{Z}_1} := \limsup_{s \to 1^+} (s-1)\nm{T}_s,
\]
where $\nm{T}_s = \Tr(|T|^s)^\frac{1}{s}$.  Then  \cite[Theorem 4.5]{CRSS}  identif\/ied that
\[
\mathcal{Z}_1 = \{ T \in \mathcal{L}^\infty | \nm{T}_{\mathcal{Z}_1}
< \infty \} \equiv \mathcal{M}_{1,\infty}
\]
and that
\[
e^{-1} \nm{T}_0 \leq \nm{T}_{\mathcal{Z}_1} \leq \nm{T}_{1,\infty},
\]
where $\nm{\cdot}_0$ is the Riesz seminorm~\eqref{eq:norm0}.

Hence, $(s-1)\zeta_{|T|}(s) \in L^\infty((1,2])$ if and only if
$T \in \mathcal{M}_{1,\infty}$, which we rewrite as (def\/ining
the following function to be $0$ for $r \in [0,1)$)
\begin{equation} \label{eq:zeta_seq0}
\frac{1}{r}\zeta_{|T|}\left(1+\frac{1}{r}\right) \in L^\infty([0,\infty))
\ \Leftrightarrow \ T \in \mathcal{M}_{1,\infty} .
\end{equation}
This is known, from \cite{LS2}, to be equivalent to
\begin{equation} \label{eq:zeta_seq}
\frac{1}{k}\zeta_{|T|}\left(1+\frac{1}{k}\right) \in \ell^\infty
\ \Leftrightarrow \ T \in \mathcal{M}_{1,\infty} .
\end{equation}
We can obtain positive unitarily invariant singular functionals
on $\mathcal{M}_{1,\infty}$, that will equate to~\eqref{eq:res0} when $T$ is a Tauberian operator, by applying a generalised limit
$\xi \in S_\infty(\ell^\infty)$ (resp.~$\phi \in S_\infty([0, \infty))$)
to the sequence \eqref{eq:zeta_seq} (resp.~function~\eqref{eq:zeta_seq0}).
Exactly which generalised limits produce \emph{linear} functionals
(and hence singular traces) is an open question.
%We now know, from a private communication by A.~Sedaev,
%that the functional
%\begin{equation} \label{eq:zeta_tr}
%\xi \left( \left\{ \frac{1}{k}\zeta_{T}(1+\frac{1}{k}) \right\}_{k=1}^\infty \right) \ , \ 0 < T \in %\mathcal{M}_{1,\infty}
%\end{equation}
%def\/ines a Dixmier trace (an element of $\mathcal{V}_1$)
%for any $\xi \in BL$.  However, we only
%know exactly which Dixmier traces in the case
%$\omega \in BL \cap DL$.
We know that choosing $\xi \in BL \cap DL$ (resp.~$\phi \in BL[0,\infty) \cap DL[0,\infty)$)
results in a Dixmier trace.

We summarise the results of \cite{LS2}, based on  \cite[Theorem~4.11]{CRSS}, see \cite{CPS} and \cite{AF} for additional information.
Def\/ine the averaging sequence $E : L^\infty([0,\infty)) \to \ell^\infty$ by
\[
E_k(f) :=  \int_{k-1}^k f(t)dt   , \qquad f \in L^\infty([0,\infty)) .
\]
Def\/ine the map $L^{-1} : L^\infty([1,\infty)) \to L^\infty([0,\infty))$ by
\[
L^{-1}(g)(t) = g(e^t)   , \qquad g \in L^\infty([1,\infty)).
\]
Def\/ine the f\/loor map
$
p : \ell^\infty \to L^\infty([1,\infty))
$
by
\[
p( \{a_k \}_{k=1}^\infty)(t) :=
\sum_{k=1}^\infty  a_k \chi_{[k,k+1)}(t)   , \qquad \{ a_k \}_{k=1}^\infty \in \ell^\infty .
\]
%$$
%p( \{a_k \}_{k=1}^\infty)(t) :=
%\sum_{k=1}^\infty  \Big( a_{k} + (a_{k+1}-a_{k})(t-k) \Big)
%\chi_{[k,k+1)}(t) .
%$$
Def\/ine, f\/inally, the mapping $\mathcal{L} : (\ell^{\infty})^*
\to (\ell^{\infty})^*$ by
\[
\mathcal{L}(\omega) := \omega \circ E \circ L^{-1} \circ p   , \qquad
\omega \in (\ell^\infty)^*.
\]
In the following statements we use the notation
established in \eqref{eq:def_BL}--\eqref{eq:def_DL} and \eqref{eq:def_BLi}--\eqref{eq:def_MDLi}.

\begin{theorem}[\protect{\cite[Theorems~3.1, 3.3]{LS2}}] \label{thm:resPA}
Let $P^* = P = P^2 \in \mathcal{L}(H)$ be a projection
and $0 < T \in \mathcal{M}_{1,\infty}$.

Then, for any $\phi \in BL[0,\infty) \cap DL[0,\infty)$,
\[
\tr_{L(\phi)}(PTP) = \phi \left( \frac{1}{r} \Tr\big(PT^{1+\frac{1}{r}}P\big) \right)  .
\]

Similarly, for any $\xi \in BL \cap DL$, $\mathcal{L}(\xi) \in D_2$ and
\[
\Tr_{\mathcal{L}(\xi)}(PTP) = \xi \left( \frac 1k \Tr\big(PT^{1+\frac 1k}P\big) \right) .
\]

Moreover, $\lim\limits_{s \to 1^+} (s-1) \Tr(PT^{s}P)$ exists iff $PTP$ is Tauberian and in either case
\[
\tr_{\upsilon}(PTP) = \Tr_{\omega}(PTP) = \lim_{s \to 1^+} (s-1) \Tr(PT^{s}P)
\]
for all $\upsilon \in DL[1,\infty)$, $\omega \in D_2$.
\end{theorem}

\begin{corollary}[\protect{\cite[Corollaries~3.2, 3.4]{LS2}}] \label{cor:resA}
Let $A \in \mathcal{L}(H)$ and $0 < T \in \mathcal{M}_{1,\infty}$.
Then, for any $\phi \in BL[0,\infty) \cap DL[0,\infty)$,
\[
\tr_{L(\phi)}(AT) = \phi \left( \frac 1r \Tr\big(A T^{1+\frac{1}{r}}\big) \right) .
\]

Similarly, for any $\xi \in BL \cap DL$,
\[
\Tr_{\mathcal{L}(\xi)}(AT) = \xi \left( \frac 1k \Tr\big(AT^{1+\frac 1k}\big) \right) .
\]

Moreover, if $PTP$ is Tauberian for all projections $P$ in the von Neumann algebra generated by $A$ and $A^*$,
\[
\tr_{\upsilon}(AT)=\Tr_{\omega}(AT) = \lim_{s \to 1^+} (s-1) \Tr(AT^s)
\]
for all $\upsilon \in DL[1,\infty)$, $\omega \in D_2$.
\end{corollary}

For the situation $s_0=p$, see \cite{CRSS}.

\subsection{Heat kernel asymptotics} \label{sec:res.2}

We follow the exposition of \cite{Sed2009}, following \cite{CPS} and \cite{CRSS}.
From \eqref{eq:def_C} and \eqref{eq:def_L_maps}
\[
L(C)(f)(t) = \frac{1}{\log t} \int_1^t f(s) \frac{ds}{s}   , \qquad f \in L^\infty([1,\infty)) .
\]
It was noted in  \cite[Lemma~5.1]{CRSS} that the function
\[
g_{T,\alpha}(t) = \frac{1}{t}\Tr( e^{-(tT)^{-\alpha}})    , \qquad \alpha > 0
\]
belongs to $L^\infty([1,\infty))$ only for $0 < T \in \mathcal{L}^{1,w}$,
where $\mathcal{L}^{1,w}$ is the symmetric subspace from \eqref{def:weak_l1} or \eqref{def:weak_lp}.  In general we have only
\[
L(C)(g_{T,\alpha}) \in L^{\infty}([1,\infty))   , \qquad
0 < T \in \mathcal{M}_{1,\infty} .
\]
For this reason heat kernel results are likely to be
restricted only to Connes--Dixmier traces ($\mathcal{V}_2$) or Ces\`{a}ro invariant Dixmier traces ($\mathcal{V}_3$)
outside of Tauberian operators in $\mathcal{L}^{1,w}$.
From~\cite{Sed2009} and~\cite{ssz2010}:

\begin{theorem}[\protect{\cite[Theorems~3--6]{Sed2009} and  \cite[Theorem~33]{ssz2010}}]
Let $\omega \in CDL[1,\infty)$ and $0 < T \in \mathcal{M}_{1,\infty}$ have trivial kernel.  Then
\[
\tr_\omega(T)
= \frac{1}{\Gamma(\frac{1}{\alpha} + 1)}
\omega \circ L(C) \left( \frac{1}{t} \Tr\big( e^{-(tT)^{-\alpha}}\big) \right)   , \qquad \alpha > 0 .
\]
Let $\omega \in CDL[1,\infty)$ and $0 < T \in \mathcal{L}^{1,w}$
have trivial kernel.
Then
\[
\tr_\omega(T)
= \frac{1}{\Gamma(\frac{1}{\alpha}+1)}
\omega \left( \frac{1}{t} \Tr\big( e^{-(tT)^{-\alpha}}\big) \right)   , \qquad \alpha > 0 .
\]
Let $\omega \in D(C)[1,\infty)$ and $0 < T \in \mathcal{M}_{1,\infty}$ have trivial kernel.
Then
\[
\tr_\omega(T)
= \frac{1}{\Gamma(\frac{1}{\alpha}+1)}
\omega \left( \frac{1}{t} \Tr\big( e^{-(tT)^{-\alpha}}\big) \right)   , \qquad \alpha > 0 .
\]
Finally, if $0 < T \in \mathcal{L}^{1,w}$ is Tauberian
with trivial kernel, then
\[
\tr_\omega(T)
= \frac{1}{\Gamma(\frac{1}{\alpha}+1)} \lim_{t \to \infty} \frac{1}{t} \Tr\big( e^{-(tT)^{-\alpha}}\big)
= \frac{1}{\Gamma(\frac{1}{\alpha}+1)} \lim_{t \to 0^+} t^{\frac{1}{\alpha}} \Tr\big( e^{-tT^{-\alpha}}\big)
  , \qquad \alpha > 0
\]
for all $\omega \in DL[1,\infty)$.
\end{theorem}

For example, if $\ds^{-n}$ is known to be Tauberian
and belongs to $\mathcal{L}^{1,w}$ then
\[
\tr_\omega(\ds^{-n})
= \frac{1}{\Gamma(\frac{n}{2}+1)} \lim_{t \to 0^+} t^{\frac{n}{2}} \Tr\big( e^{-tD^2}\big) .
\]
From Weyl's formula for the eigenvalues of the Hodge Laplacian on a $n$-dimensional compact Riemannian manifold (and also using comments from  \cite[p.~278]{CRSS}), we have that
\[
0 < (1+Q^2)^{-n/2d} \in \mathcal{L}^{1,w}
\]
and is Tauberian for any positive elliptic operator in $Q \in \mathrm{Op}_{\mathrm{cl}}^d(M)$.  This fact, combined with the above
equations, reconstructs~\eqref{eq:htk_classic}.

For heat kernel formulas involving $\tr_\omega(AT)$, $A \in \mathcal{L}(H)$,
$0 < T \in \mathcal{M}_{1,\infty}$,
see \cite{CPS} and \cite{CRSS}.

\section{Characterising the noncommutative integral} \label{sec:char}

The functional
\begin{equation} \label{eq:noncomm_meas_g}
\Phi_{\mathrm{f},T}(a) := \mathrm{f}(aT)   , \qquad a \in \mathcal{L}(H),
\end{equation}
where $\mathrm{f}$ is a trace on a two-sided ideal of compact operators $\mathcal{J}$ and $0 < T \in \mathcal{J}$,
is the general format of noncommutative integration
introduced by Connes.
We are concerned in this section
with the characterisation of
linear functionals on $\mathcal{L}(H)$ constructed by~\eqref{eq:noncomm_meas_g}.

The functional $\Phi_{\mathrm{f},T}$ is normal (contained in
$\mathcal{L}(H)_*$) if and only if $\Phi_{\mathrm{f},T}(a) = \Tr(aT)$,
$\fa\, a \in \mathcal{L}(H)$, for a positive trace class `density' $T$.
From the def\/inition of the canonical trace $\Tr$:
\[
\Tr(aT)  = \sum_{m=1}^\infty \lambda_m \inprod{h_m}{ah_m},
\]
where $\{ h_m \}_{m=1}^\infty$ is an orthonormal basis of eigenvectors for the positive compact operator $T$ with eigenvalues $0 < \{ \lambda_m \}_{m=1}^\infty \in \ell^1$.
Def\/ine the normal linear functional on $\ell^\infty$,
$\sigma_T \in (\ell^\infty)_* \cong \ell^1$, by
\[
\sigma_T(\{a_m\}_{m=1}^\infty) := \sum_{m=1}^\infty \lambda_m a_m .
\]
Then
\begin{equation} \label{eq:tr_char}
\Tr(aT) = \sigma_T(\{\inprod{h_m}{ah_m}\}_{m=1}^\infty)   , \qquad \fa\, a \in \mathcal{L}(H).
\end{equation}
Conversely,
\[
\sigma(\{ \inprod{h_m}{ah_m} \}_{m=1}^\infty)   , \qquad 0 < \sigma \in (\ell^\infty)_*   , \qquad a \in \mathcal{L}(H)
\]
def\/ines a positive normal linear functional on $\mathcal{L}(H)$ for any choice of orthonormal basis $\{ h_m \}_{m=1}^\infty$ of $H$.
If $\mathbf{e}_m$ is the sequence with $1$ in the $m^\mathrm{th}$ place and $0$ otherwise, set $T_\sigma$
to be the operator with eigenvalues $0 < \{ \sigma(\mathbf{e}_m) \}_{m=1}^\infty \in \ell_1$ associated to the basis $\{ h_m \}_{m=1}^\infty$.  By this construction
\begin{equation} \label{eq:tr_char_2}
\Tr(aT_\sigma) = \sigma(\{\inprod{h_m}{ah_m}\}_{m=1}^\infty)   , \qquad \fa \, a \in \mathcal{L}(H).
\end{equation}

By \eqref{eq:tr_char} and \eqref{eq:tr_char_2} the positive normal linear functionals
\[
\Phi_{\Tr,T}(a) = \Tr(aT)   , \qquad 0 < T \in \mathcal{L}^1
\]
are exactly characterised by
\[
0 < \sigma \in (\ell^\infty)_* \ \ \mbox{applied to the sequence of ``expectation values''}  \ \
\{ \inprod{h_m}{ah_m} \}_{m=1}^\infty
\]
for some orthonormal basis $\{h_m\}_{m=1}^\infty$.
We now know that a similar characterisation exists for any functional
$\Phi_{\mathrm{f},T}$ as def\/ined in~\eqref{eq:noncomm_meas_g}.

\subsection{Characterisation and singular traces}

Let $\mathcal{J}$ be a (two sided) ideal contained in the compact operators $\mathcal{L}^\infty$ of the separable complex Hilbert space $H$.  By a
% continuous
trace on $\mathcal{J}$ we mean a linear functional $\mathrm{f} : \mathcal{J} \to \CC$ such that
% $|\mathrm{f}(aT)| \leq \nm{a}|\mathrm{f}(T)|$ and
$\mathrm{f}([a,T])=0$ for all $a\in \mathcal{L}(H)$, $T \in \mathcal{J}$, i.e.~$\mathrm{f}$ vanishes on the commutator subspace
$\mathrm{Com} \, \mathcal{J}$ of $\mathcal{J}$.
Note it is not assumed that $\mathrm{f}$ is positive or continuous.

The characterisation for \eqref{eq:noncomm_meas_g} in
Theorem~\ref{thm:char} below
was f\/irst shown for
$\mathrm{f} \in \mathcal{V}_{PDL[1,\infty)}$, $0 < T \in \mathcal{M}_{1,\infty}$, in \cite[Theorem~3.6]{LS2} using residues (Corollary~\ref{cor:resA}).  The result below (to appear)
is due to collaboration with N.~Kalton and uses results on the commutator subspace and sums of commutators
\cite{K1998, DFWW2004}.

\begin{theorem} \label{thm:char}
Let $0 < T \in \mathcal{J}$ and $\mathrm{f}$ be a trace on $\mathcal{J}$.  Then there exists an
orthonormal basis $\{ h_m \}_{m=1}^\infty$
of $H$ $($consisting of eigenvectors for $T)$ and a linear functional
$\sigma : \ell^\infty \to \CC$ such that
\begin{equation} \label{eq:char}
\mathrm{f}(aT) = \sigma( \{ \inprod{h_m}{ah_m} \}_{m=1}^\infty )
\end{equation}
for all $a \in \mathcal{L}(H)$.
\end{theorem}

Straightforward corollaries show, when $\mathrm{f} > 0$,
that $\mathrm{f}$ has the property
\[
|\mathrm{f}(aT)| \leq C \nm{a} \ \ \mbox{some} \ \ C > 0,
\]
if and only if $0 < \sigma \in (\ell^\infty)^*$.
The basis $\{ h_m \}_{m=1}^\infty$
in the theorem is any basis of eigenvectors of~$T$ rearranged so
that the sequence of eigenvalues $\lambda_m$ associated
to $h_m$ is a decreasing sequence.

A singular trace is a trace on a two sided ideal $\mathcal{J}$, containing the f\/inite rank operators $F$, that vanishes on $F \subset \mathcal{J}$.
The theorem has the following corollary.

\begin{corollary} \label{cor:char}
Let $0 < T \in \mathcal{J}$ and $\mathrm{f}$ be a positive
singular trace on $\mathcal{J}$.
Then there exists an orthonormal basis $\{ h_m\}_{m=1}^\infty$
of $H$ $($consisting of eigenvectors for $T)$ and a singular state
$($generalised limit$)$ $L : \ell^\infty \to \CC$ such that
\begin{equation} \label{eq:char_sing}
\mathrm{f}(aT) = \mathrm{f}(T)L( \{ \inprod{h_m}{ah_m} \}_{m=1}^\infty )
\end{equation}
for all $a \in \mathcal{L}(H)$.
\end{corollary}

In particular, for the noncommutative integral from \eqref{eq:ncg_state_normal},
\[
\Phi_\rho(a) = \rho(a\ds^{-n})  , \qquad a \in \A'', \qquad \rho \in \mathcal{V}_{\mathrm{sing}}^+
  , \qquad \rho(\ds^{-n}) = 1
\]
the corollary implies that there exists a basis of eigenvectors
$\{ h_m\}_{m=1}^\infty$ of $D$ such that
\[
\rho(a\ds^{-n}) = L_\rho( \{ \inprod{h_m}{ah_m} \}_{m=1}^\infty )
\]
where $L_\rho$ is a generalised limit.

\begin{example} \label{ex:torus_integral}
Let $\TT^n$ be the $n$-torus, $\mathbb{R}^n / \mathbb{Z}^n$, with Laplacian $\Delta$.
Let $\mathcal{J}$ be any two-sided ideal with
trace $\mathrm{f}$ and $g$ be a bounded Borel function $g$
such that $g(\Delta) \in \mathcal{J}$.
Let $f \in L^\infty(\TT^n)$
and denote by
\[
f g(\Delta)
\]
the action of $f$ on $L^2(\TT^n)$ by pointwise multiplication
coupled with the compact operator $g(\Delta)$, c.f.~\cite[\S~4]{S}. Then, by Theorem~\ref{thm:char},
\[
\mathrm{f}( f g(\Delta) )
= \sigma( \inprod{e^{i \mathbf{m} \cdot x}}{fe^{i \mathbf{m} \cdot x}} )
\]
where $\mathbf{m} \in \ZZ^n$ (ordered by Cantor enumeration)
and $\sigma$ is a linear functional on $\ell^\infty$.
As
\[
\inprod{e^{i \mathbf{m} \cdot x}}{f e^{i \mathbf{m} \cdot x}}
= \int_{\TT^n} e^{-i\mathbf{m} \cdot x} f(x) e^{i\mathbf{m} \cdot x} dx
= \int_{\TT^n} f(x) dx
\]
and
\[
\sigma( \mathbf{1} ) = \mathrm{f}(g(\Delta))
\]
we obtain
\begin{equation} \label{eq:torus_int}
\mathrm{f}( f g(\Delta) ) = \mathrm{f}(g(\Delta) ) \int_{\TT^n} f(x) dx .
\end{equation}
Therefore, for any ideal $\mathcal{J}$ with a non-trivial trace $\mathrm{f}$,
by choosing the appropriate bounded Borel function $g$,
the expression
\[
\mathrm{f}( \cdot g(\Delta) )
\]
may serve as the integral on the $n$-torus (up to a constant).
Note the identif\/ication \eqref{eq:torus_int}
holds for all $f \in L^\infty(\TT^n)$.
\end{example}

Can the same statement as for the $n$-torus in the example
above be made for the Lebesgue integral of a compact Riemannian manifold $M$,  considering non-trivial positive singular traces instead of arbitrary traces?

We certainly know, from Corollary \ref{cor:char} and Theorem \ref{thm:result2}, that if the sequence
\begin{equation} \label{eq:part_QUE}
\{ \inprod{h_m}{fh_m} \}_{m=1}^\infty \in c   , \qquad \fa \, f \in C^\infty(M)
\end{equation}
for a basis of eigenvectors $\{ h_m\}_{m=1}^\infty$ of the Hodge--Laplacian $\Delta$ (ordered so that the eigenvalues associated to the eigenvectors are decreasing), then
\[
\mathrm{f}( f g(\Delta) ) = \mathrm{f}(g(\Delta))
\int_M f dx   , \qquad \fa f \in C^\infty(M)
\]
for any positive non-trivial singular trace $\mathrm{f}$ on
an ideal $\mathcal{J}$ and with an
appropriate choice of bounded Borel function $g$
(so that $g(\Delta) \in \mathcal{J}$).

The property \eqref{eq:part_QUE} is a
restricted form of a property for compact Riemannian manifolds
known as Quantum Unique Ergodicity (QUE)
\cite{Donn03, Sarnak1}.  It is not well known which manifolds have the QUE property.

The same result as for the $n$-torus applies to anything `f\/lat'.
The next example involving the non-commutative torus f\/irst appeared
(for $\mathrm{f} \in \mathcal{V}_{PDL[1,\infty)}$) in \cite{LS2}.

\begin{example}{\sloppy
Consider two unitaries $u,v$ such that $uv = \lambda vu$,
for $\lambda := e^{2\pi i \theta} \in \SB$ (the unit circle). Denote by
$F_\theta(u,v)$ the $*$-algebra of linear combinations
$\sum_{(m,n) \in J} a_{m,n} u^mv^n$, $J \subset \ZZ^2$ is a f\/inite set, with product $ab = \sum_{r,s}( \sum_{m,n} a_{r-m,n} \lambda^{mn} b_{m,s-n} ) u^rv^s$
and involution $a^* = \sum_{r,s} (\lambda^{rs} \overline{a}_{-r,-s}) u^rv^s$, $a,b \in F_\theta(u,v)$. The assignment $\tau_0(a) = a_{0,0}$ is a faithful trace on $F_\theta(u,v)$.  Let $(H_\theta,\pi_\theta)$ denote the cyclic representation associated to $\tau_0$.
The closure, $C_\theta(u,v)$, of $\pi_\theta(F_\theta(u,v))$ in the operator norm is called a rotation $C^*$-algebra, \cite{Rieff1},
or the noncommutative torus ($\lambda \not=1$) \cite{CN, C4, GBVF}.   Canonically, f\/inite linear combinations of
$u^mv^n \hookrightarrow H_\theta$ are dense in $H_\theta$.
Def\/ine $\Delta_\theta(u^mv^n) = (m^2 + n^2)u^mv^n$.
It can be shown that the `noncommutative Laplacian' $\Delta_\theta$ has a unique positive extension (also denoted $\Delta_\theta$)
$\Delta_\theta : \Dom(\Delta_\theta) \to H_\theta$
with compact resolvent, see the previous citations.  The eigenvectors $h_{m,n} = u^mv^n \in H_\theta$ form a complete orthonormal system.
Note that
\[
\inprod{h_{m,n}}{\pi_\theta(a)h_{m,n}} =:  \tau_0((u^mv^n)^*au^mv^n)
= \tau_0(a)   , \qquad \fa \, a \in F_\theta(u,v)
\]
for any $(m,n) \in \ZZ^2$.
Using the Cantor enumeration of $\ZZ^2$, from Theorem~\ref{thm:char} we obtain
\[
\mathrm{f}(ag(\Delta_\theta))
=  \mathrm{f}(g(\Delta_\theta)) \tau_0(a)   , \qquad \fa \, a \in F_\theta(u,v)
\]
for any ideal $\mathcal{J}$ with a non-trivial trace $\mathrm{f}$, by choosing the appropriate bounded Borel function~$g$
such that $g(\Delta_\theta) \in \mathcal{J}$.
In fact,
$\inprod{h_{m,n}}{b h_{m,n}} = \tau(b)$
for any $b$ in the weak closure $C_\theta(u,v)''$
of $\pi_\theta(F_\theta(u,v))$ and where $\tau$ denotes the
normal extension of $\tau_0$.  Hence
\[
\mathrm{f}(ag(\Delta_\theta))
=  \mathrm{f}(g(\Delta_\theta)) \tau(a)   , \qquad \fa\, a \in C_\theta(u,v)'' .
\]}
\end{example}

\subsection{A dominated convergence theorem for the noncommutative integral} \label{sec:char:dct}

The following dominated convergence theorem
appeared in~\cite{LS2} as a statement for Dixmier traces.
We provide the general statement here.

Let $\mathcal{N}$ be a weakly closed $*$-subalgebra of $B(H)$.
We recall that $\mathcal{N}_*$ denotes the predual of~$\mathcal{N}$, or the set of all normal linear functionals on~$\mathcal{N}$.

We say a positive compact operator
$T$ is $(\mathcal{N},h)$-dominated if, for some orthonormal basis $\{h_m\}_{m=1}^\infty$ of eigenvectors of $T$,
there exists $h \in H$ such that $\nm{Ph_m} \leq \nm{Ph}$ for all projections $P \in \mathcal{N}$.

\begin{theorem} \label{thm:normM}
Let $0 < T \in \mathcal{J}$ and $\mathrm{f}$ be a positive singular trace on $\mathcal{J}$.  If $T$ is $(\mathcal{N},h)$-dominated,
then $\Phi_{\mathrm{f},T} = f( \cdot T) \in \mathcal{N}_*$.
\end{theorem}
\begin{proof}
Without loss $\{ h_m \}_{m=1}^\infty$ can be rearranged so that
the eigenvalues $\lambda_m$ associated to $h_m$ form a decreasing
sequence.
By hypothesis $\inprod{h_m}{Ph_m} \leq \inprod{h}{Ph}$ for all projections $P \in \mathcal{N}$.  Then $\inprod{h_m}{ah_m} \leq \inprod{h}{ah}$, $0 < a \in \mathcal{N}$, as $a$
is a uniform limit of f\/inite linear positive spans of projections
 \cite[p.~23]{Ped}.  For any generalised limit
$L$ and $0 < a \in \mathcal{N}$,
\begin{equation} \label{eq:dom}
L(\inprod{h_m}{ah_m})
\leq \limsup_{m \to \infty} \inprod{h_m}{ah_m}
\leq \inprod{h}{ah}.
\end{equation}
Let $\{ a_\alpha \}$ be a net of monotonically increasing positive elements of $\mathcal{N}$ with upper bound.  It follows that $\{ a_\alpha \}$ converges strongly
to a l.u.b.~$a \in \mathcal{N}$ \cite[p.~22]{Ped}.  From \eqref{eq:dom}
$
L(\inprod{h_m}{(a - a_\alpha)h_m})
\leq \inprod{h}{(a-a_\alpha)h}.
$
Since $\inprod{h}{(a-a_\alpha)h} \stackrel{\alpha}{\to} 0$,
$L(\inprod{h_m}{ah_m}) = \sup_{\alpha} L(\inprod{h_m}{a_\alpha h_m})$.  From Corollary~\ref{cor:char}
$\Phi_{\mathrm{f},T}(a) = \sup_\alpha \Phi_{\mathrm{f},T}(a_\alpha)$
and $\Phi_{\mathrm{f},T}$ is normal on $\mathcal{N}$  \cite[\S~3.6.1]{Ped}.
\end{proof}

We close with an example of $(\mathcal{N},h)$-dominated compact operators.

Let $H = L^2(X,\mu)$ for a $\sigma$-f\/inite
measure space $(X,\mu)$ and set $\mathcal{N}$ as multiplication operators of $L^\infty$-functions.  Then $0 < T \in \mathcal{L}^\infty(L^2(X,\mu))$ being $(\mathcal{N},h)$-dominated
is the same as $\int_J |h_m(x)|^2 d\mu(x)
\leq \int_J |h(x)|^2 d\mu(x)$ for all measurable sets $J \subset X$,
which is equivalent to the statement $|h_m|^2$ are
dominated by some $|h|^2 \in L^1(X,\mu)$ $\mu$-a.e.,
where $\{ h_m \}_{m=1}^\infty$ is any basis of eigenvectors
of $T$.  Hence any positive compact operator
$0 < T \in \mathcal{L}^\infty(L^2(X,\mu))$ is
$(\mathcal{N},\sqrt{g})$-dominated if it has a basis
of eigenvectors whose modulus squared are dominated by some positive integrable function
$g \in L^1(X,\mu)$.  For instance, the negative powers
of the Laplacian on the $n$-torus are dominated by the constant function $1 \in L^1(\TT^n)$.

\section{Summary}

We summarise the main points of the review and list the open
questions raised.

\subsection{Summary}

There are three main non-identical classes of Dixmier traces.
They are, in descending order of inclusion:
the (original) Dixmier traces ($\mathcal{V}_1$);
the Connes--Dixmier traces ($\mathcal{V}_2$); and
the Ces\`{a}ro invariant Dixmier traces ($\mathcal{V}_3$).

The notion of measurable operator with respect to each of these
classes are non-identical.  The notion of Tauberian, \eqref{eq:Taub},
is the strongest notion with ascending order of inclusion
\[
\mathcal{T} \subsetneq \mathcal{K}_1
\subsetneq \mathcal{K}_2
\subsetneq \mathcal{K}_3,
\]
where
\[
\mathcal{K}_i = \{ T \in \mathcal{M}_{1,\infty} \,|\,
\mathrm{f}(T) = \mathrm{const} \ \fa \, \mathrm{f} \in \mathcal{V}_i \}
  , \qquad i=1,2,3 .
\]
None of these sets is an ideal of compact operators.
The notion of measurable operator using Dixmier traces does not achieve $\mathcal{T} = \mathcal{K}_1$, however $\mathcal{T}^+ = \mathcal{K}_1^+$.
This has the consequence that the operators $T \in \mathcal{K}_1$ do not decompose into unique measurable components.
In short, $\mathcal{K}_1$ does not have the property
$T \in \mathcal{K}_1 \Rightarrow |T| \in \mathcal{K}_1$.

Connes' trace theorem (Theorem \ref{thm:2}) is valid for any
$\rho \in \mathcal{V}_{\mathrm{sing}}^+$, \eqref{def:V_sing}.  In particular, if $M$ is an $n$-dimensional compact Riemannian manifold with Hodge Laplacian $\Delta$ we have
\[
\rho(P) = c_\rho \mathrm{Res}(P) \qquad \fa \, P \in \mathrm{Op}^{-n}_{\mathrm{cl}}(M)
\]
with the constant $c_\rho>0$ for the subset $\Theta_n$ of $\mathcal{V}_{\mathrm{sing}}^+$ such that $\rho((1+\Delta)^{-n/2}) > 0$.  The set $\Theta_n$ of singular traces is larger than the set of Dixmier traces ($\mathcal{V}_1$).
Further, for $\rho \in \Theta_n$ we have
\[
\rho(f(1+\Delta)^{-n/2})
= \frac{\rho((1+\Delta)^{-n/2})}{\mathrm{Vol}(M)} \int_{M} f dx   , \qquad \fa f \in L^\infty(M),
\]
where $dx$ is the volume form on $M$ (Theorem \ref{thm:3}).

For any spectral triple $(\A,H,D)$ with $\ds^{-n} \in \mathcal{M}_{1,\infty}$
def\/ine
$\Theta(D)$ to be the subset of $\mathcal{V}_{\mathrm{sing}}^+$ such that $\rho(\ds^{-n}) > 0$.  Def\/ine the noncommutative integral
\[
\rho(a\ds^{-n})   , \qquad a \in \A   , \qquad \rho \in \Theta(D) .
\]
This def\/ines a family of continuous positive traces on $\A$
(in the implied $C^*$-norm) (Corollary~\ref{cor:trace}).

Restricting the noncommutative integral to the consideration
of $\rho \in \mathcal{V}_1$ (the set of Dixmier traces
or the set of Hardy--Littlewood--P\'{o}lya submajorisation
ordered continuous functionals on~$m_{1,\infty}$, Section~\ref{sec:dtr:ssf}) enables Lidksii type theorems,
residue, and heat kernel theorems, see Section~\ref{sec:res}
and comments.

The functional
\[
\mathrm{f}(aT)   , \qquad a \in \A \subset \mathcal{L}(H)
\]
where $\mathrm{f}$ is a trace on a two-sided ideal of compact operators $\mathcal{J}$ and $0 < T \in \mathcal{J}$,
is the general format of noncommutative integration
introduced by Connes.  Such functionals can be characterised by
linear functionals $\sigma : \ell^\infty \to \mathbb{C}$
applied to sequences of ``expectation values''
$\{ \inprod{h_m}{ah_m} \}_{m=1}^\infty$ for some orthonormal basis $\{h_m\}_{m=1}^\infty$ of eigenvectors of $T$.  In particular
\[
\mathrm{f}(aT) = \sigma( \{ \inprod{h_m}{ah_m} \}_{m=1}^\infty )
\]
for all $a \in \mathcal{L}(H)$.

For f\/lat torii, both commutative and noncommutative, this implies
\[
\mathrm{f}(ag(\Delta)) = \mathrm{f}(g(\Delta)) \tau(a)
\]
for any trace $\mathrm{f}$ on a two-sided ideal $\mathcal{J}$,
and $g$ a positive bounded Borel function such that
$g(\Delta) \in \mathcal{J}$.  Here $\Delta$ is the
Laplacian (resp.~noncommutative Laplacian)
on the torus and $\tau$ is the Lebesgue integral on commutative torii (resp.~the unique faithful normal f\/inite trace on the type II$_1$
noncommutative torii).

If $\mathcal{N} \subset \mathcal{L}(H)$ is a von Neumann algebra,
a positive compact operator
$T$ is $(\mathcal{N},h)$-dominated if, for some orthonormal basis $\{h_m\}_{m=1}^\infty$ of eigenvectors of $T$,
there exists $h \in H$ such that $\nm{Ph_m} \leq \nm{Ph}$ for all projections $P \in \mathcal{N}$.
If $\mathcal{J}$ is any two-sided ideal with positive singular trace
$\mathrm{f}$, and $0 < T \in \mathcal{J}$ is $(\mathcal{N},h)$-dominated,
then $\mathrm{f}(\cdot T)$ is a positive normal linear functional
on~$\mathcal{N}$.

\subsection{List of open questions}

Let $(\A,H,D)$ be a spectral triple (as def\/ined in the f\/irst paragraph of Section~\ref{sec:intro}) such that $\ds^{-n} \in \mathcal{M}_{1,\infty}$ (as def\/ined at \eqref{eq:intro_dix_ideal}).

\begin{openq}[from Section~\ref{sec:intro}] If $\A''$ contains no finite rank operators is $\Phi_\omega$ $($the functional in \eqref{eq:ncg_int}$)$ normal?
\end{openq}

\begin{openq}[from Section~\ref{sec:intro}]  What are necessary and sufficient relationships between $\ds^{-n}$ and the $*$-algebra $\A$
so that $\Phi_\omega$ $($the functional in \eqref{eq:ncg_int}$)$
is independent of $\omega$ $($for various sets of generalised limits$)$?
\end{openq}

\begin{openq}[from Section~\ref{sec:meas}] Is $B \subset S_\infty(\ell^\infty)$
satisfying the Ces\`{a}ro limit property $($from~{\rm \cite[\S~5]{LSS})}
a necessary and sufficient condition for equality of a set
$\mathcal{K}_B^+$ of positive measurable operators $($e.g.~\eqref{eq:K}$)$ with the Tauberian operators?
\end{openq}

\begin{openq}[from Section~\ref{sec:meas}] Is $\mathcal{K}_{PDL[1,\infty)}^+$ $($defined in~\eqref{eq:K}$)$ equal to the Tauberian operators?
\end{openq}

\begin{openq}[from Section~\ref{sec:dtr:subideals}] Is there a maximal closed symmetric subideal $\mathcal{U}_2$ of~$\mathcal{M}_{1,\infty}$ such that $\tilde{\mathcal{T}} \cap \mathcal{U}_{2} = \mathcal{K}_2 \cap \mathcal{U}_{2}$?  Is $\mathcal{L}^{1,w} \subseteq \mathcal{U}_2$?
\end{openq}

%(Sect.~\ref{sec:res.1}) Which generalised limits produce \emph{linear} functionals
%(and hence singular traces) on $\mathcal{M}_{1,\infty}$
%from the formula $\rho_\xi$ (defined in \eqref{eq:zeta_tr})?  What is the relation of this set of %singular traces to $\mathcal{V}_{\frac{1}{2}}$ (defined in~\eqref{eq:dtr0.5})?

\begin{openq} What is the relation between the sets $\mathcal{V}_{\frac{1}{2}}$ and $\Theta(D)$ $($defined above$)$ of singular traces?
What are their sets of measurable operators?
\end{openq}

\appendix

\section{Identif\/ications in (\ref{fig:simplify})} \label{sec:app_A}

We start with preliminaries.  Def\/ine
\begin{equation} \label{eq:app:alpha}
\alpha_g(t) := \frac{1}{\log(1+t)} \int_1^t g^*(s)ds   , \qquad t \geq 1,
\end{equation}
where
\begin{equation} \label{eq:app:g*}
g^*(s) := \inf \{ t | \mu(|g| > t) < s \}   , \qquad \mu \ \mbox{is Lebesgue measure on $[1,\infty)$}.
\end{equation}
Set
\begin{equation}
\mathrm{m}_{1,\infty} := \{ g \in L^\infty([1,\infty))
| \alpha_g \in C_b([1,\infty)) \} .
\end{equation}
Def\/ine: $p$ from $\ell^\infty$ to bounded Borel functions on $[0,\infty)$ by
\begin{equation} \label{eq:app:def_p}
p( \{a_k \}_{k=0}^\infty)(t) :=
\sum_{k=0}^\infty  a_k \chi_{[k,k+1)}(t)  , \qquad t \geq 0 ;
\end{equation}
$p_c$ from $\ell^\infty$ to continuous bounded functions on $[0,\infty)$ by
\begin{equation} \label{eq:app:def_pc}
p_c( \{a_k \}_{k=0}^\infty)(t) :=
\sum_{k=0}^\infty \left( a_{k} + (t-k)(a_{k+1} - a_{k}) \right) \chi_{[k,k+1)}(t)   , \qquad t \geq 0 ;
\end{equation}
$r$ from bounded Borel or continuous functions on $[0,\infty)$ to $\ell^\infty$ by
\begin{equation} \label{eq:app:def_r}
r(f)_k := f(k)   , \qquad k \in \{0,1,2,\ldots\};
\end{equation}
$E : L^\infty([0,\infty)) \to L^\infty([0,\infty))$ by
\begin{equation} \label{eq:app:E}
E(f)(t) :=  \int_{t}^{t+1} f(s)ds   , \qquad t \geq 0 ;
\end{equation}
$L^{-1} : L^\infty([1,\infty)) \to L^\infty([0,\infty))$ by
\begin{equation} \label{eq:app:L}
L^{-1}(g)(t) = g(e^t)   , \qquad t \geq 0 ;
\end{equation}
and f\/inally the maps on $\ell^\infty$ for $j \in \NN$ by
\begin{gather}
T_j( \{a_k\}_{k=0}^\infty )   :=   \{ a_{k+j} \}_{k=0}^\infty , \label{eq:app:def_Td} \\
C( \{a_k \}_{k=0}^\infty)   :=
\left\{ \frac{1}{n+1} \sum_{k=0}^n  a_k \right\}_{n=0}^\infty, \label{eq:app:def_C}
\end{gather}
and the maps on $L^\infty([0,\infty))$ for $a > 0$ by
\begin{gather}
T_a(f)(t)   :=   f(t+a), \label{eq:app:def_Tc} \\
C(f)(t)   :=   \frac{1}{t} \int_0^t f(s)ds. \label{eq:app:def_Cc}
\end{gather}

\begin{lemma} \label{lemma:app:prelim}
Let $\{a_n \}_{n=0}^\infty \in \ell^\infty$ and $f \in L^\infty([0,\infty))$.
With the above definitions we have:
\begin{enumerate}\itemsep=0pt
\item $\lim_{t \to \infty} (T_jp_{(c)}-p_{(c)}T_j)(\{a_n \})(t) = 0$;
\item $\lim_{t \to \infty} (Cp_{(c)}-p_{(c)}C)(\{a_n \})(t) = 0$;
\item $(T_jrE - rET_j)(f) \in c_0$;
\item $(CrE - rEC)(f) \in c_0$,
\end{enumerate}
for any $j \in \NN$.
\end{lemma}
\begin{proof}
Let $\{ a_k \}_{k=0}^\infty \in \ell^\infty$. Then
\[
(T_jp-pT_j)(\{ a_k \})(t) =  \sum_{k=0}^\infty  a_k \chi_{[k,k+1)}(t+j)
- \sum_{k=0}^\infty  a_{k+j} \chi_{[k,k+1)}(t)
 =  0 .
\]
Similarly,
\[
(T_jp_c-p_cT_j)(\{ a_k \})(t) = 0 .
\]
This proves 1.
We also have
\begin{gather*}
(Cp-pC)(\{ a_k \})(t)   =   \frac{1}{t} \int_0^t \sum_{k=0}^\infty  a_k \chi_{[k,k+1)}(s) ds - \sum_{k=0}^\infty  \frac{1}{k+1} \sum_{i=0}^k a_i \chi_{[k,k+1)}(t) \\
 \phantom{(Cp-pC)(\{ a_k \})(t)}{}
  =    \left(\frac{\floor{t}+1}{t}-1\right)(C (\{ a_k \})(\floor{t}) +
  \left(1-\frac{\floor{t}}{t}\right) a_{\floor{t}}.
\end{gather*}
Hence
\[
\nm{(Cp-pC)\{ a_k \}}_\infty
\leq \left(\left|\frac{\ceil{t}}{t}-1\right|+\left|\frac{\floor{t}}{t}-1\right|\right) \nm{a}_\infty \to 0 .
\]
Similarly,
\begin{gather*}
(Cp_c-p_cC)(\{ a_k \})(t)   =   \frac{1}{t}
\int_0^t \sum_{k=0}^\infty (a_{k} + (s-k)(a_{k+1}-a_{k})) \chi_{[k,k+1)}(s)ds \\
\qquad\quad{}
- \sum_{n=0}^\infty C (\{ a_k \})(n)
+(t-n)(C (\{ a_k \})(n+1) - C (\{ a_k \})(n)) \chi_{[n,n+1)}(t) \\
\qquad {} =    \frac{1}{t} \sum_{k=0}^{\floor{t}} \frac{1}{2}(a_{k} + a_{k+1}) - \frac{1}{t} \int_t^{\ceil{t}} a_{\floor{t}}
+(s-\floor{t})(a_{\floor{t}+1}-a_{\floor{t}})ds \\
 \qquad\quad{} - C (\{ a_k \})(\floor{t})
-(t-\floor{t})(C (\{ a_k \})(\floor{t}+1) - C (\{ a_k \})(\floor{t})).
\end{gather*}
We recall that
\[
C(\{a_k\})(n+1) - C(\{a_k\})(n) =
\left(\frac{n+1}{n+2}-1\right) C(\{a_k\})(n) - \frac{a_{n+1}}{n+2}
\]
so that
\[
\lim_{n \to \infty} |C(\{a_k\})(n+1) - C(\{a_k\})(n)| = 0.
\]
Also
\[
C\left(\frac{1}{2}(\{a_{k+1}+a_{k} \})\right)(n) =
\frac{1}{2}\left(\frac{a_{n+1}}{n+1} - \frac{a_{0}}{n+1}\right)
+ C(\{ a_k \})(n)
\]
so that
\[
\lim_{n \to \infty} \left|C\left(\frac{1}{2}(\{a_{k+1}+a_{k} \})\right)(n)
 - C(\{ a_k \})(n)\right| = 0.
\]
Now
\begin{gather*}
\nm{(Cp_c-p_cC)\{ a_k \}}_\infty
  \leq   \left|\frac{\ceil{t}}{t}-1\right| \nm{a}_\infty
+ \left|C\left \{\frac{1}{2}(a_{k+1}+a_k)\right\}(\floor{t})
- C (\{ a_k \})(\floor{t})\right| \\
\phantom{\nm{(Cp_c-p_cC)\{ a_k \}}_\infty   \leq }{}
 + |C (\{ a_k \})(\floor{t}+1) - C (\{ a_k \})(\floor{t})| + 3 \frac{\nm{a}_\infty}{t} \to 0 .
\end{gather*}
These results demonstrate 2.

Set $E' = rE : L^\infty([0,\infty)) \to \ell^\infty$.  Consider
\[
(T_jE' - E'T_j)f(n) = \int_{n}^{n+1} f(t+j)dt - \int_{n+j}^{n+j+1} f(t)dt = 0 .
\]
We have
\begin{gather*}
|(E'C-CE')f(n)|   =   \left|\int_{n}^{n+1} \frac{1}{t} \int_0^t f(s)ds dt - \frac{1}{n+1} \sum_{i=0}^n \int_{i}^{i+1} f(s)ds\right| \\
 \phantom{|(E'C-CE')f(n)|}{}
 =    \left|\int_{n}^{n+1} \frac{1}{t} \int_0^t f(s)ds dt - \frac{1}{n+1} \int_{0}^{n+1} f(s)ds\right| \\
 \phantom{|(E'C-CE')f(n)|}{}
   \leq    \sup_{t \in [n,n+1)}\left|\frac{1}{t} \int_0^t f(s) ds - \frac{1}{n+1} \int_{0}^{n+1} f(s)ds\right| \\
\phantom{|(E'C-CE')f(n)|}{}  \leq    \left(\frac{n+1}{n}-1+\frac{1}{n}\right)\nm{f}_\infty \to 0 .
\end{gather*}
These equations demonstrate 3 and 4.
\end{proof}

\begin{lemma} \label{lemma:app:prelim2}
If $g \in \mathrm{m}_{1,\infty}$, then:
\begin{enumerate}
\item $\xi(1-p_{(c)}r)L^{-1}(\alpha_g) = 0$ for all $\xi \in BL_{(c)}[0,\infty)$;
\item $\xi'r(1 - E)L^{-1}(\alpha_g) = 0$ for all $\xi' \in BL$.
\end{enumerate}
\end{lemma}
\begin{proof}
The following arguments were f\/irst published in a more general form in~\cite[\S~2.2]{LSS}.
We include them for completeness.
Let $f = L^{-1}(\alpha_g) \in C_b([0,\infty))$ and $c(s) = \log(1+e^s)$.  Set
\[
k(s) = \frac{1}{c(s)} \int_{e^s}^{e^{s+1}} g^*(u)du    , \qquad s \geq 0 .
\]
Let $\theta$ be a state of $C_b([0,\infty))$ such that $\theta = \theta T_j$
for $j \in \NN$.  For example, $\theta = \xi|_{C_b([0,\infty))}$ where $\xi \in BL[0,\infty)$,
$\theta = \xi$ where $\xi \in BL_c[0,\infty)$, or $\xi'r$ for $\xi' \in BL$.  Then
\begin{gather*}
\theta(k)   =   \theta \left( \frac{1}{c(s)} \int_{e^s}^{e^{s+1}} g^*(u)du \right)
  =   \theta \left( \frac{1}{c(s)} \int_{1}^{e^{s+1}} g^*(u)du \right)
- \theta \left( \frac{1}{c(s)} \int_{1}^{e^{s}} g^*(u)du \right) \\
\phantom{\theta(k)}{}
  =   \theta \left( \frac{1}{c(s)} \int_{1}^{e^{s+1}} g^*(u)du \right)
- \theta \left( \frac{1}{c(s+1)} \int_{1}^{e^{s+1}} g^*(u)du \right) \\
\phantom{\theta(k)}{}  =   \theta \left((1 - \frac{c(s+1)}{c(s)})f(s+1) \right)
  =   \lim_{s \to \infty} (1 - \frac{c(s+1)}{c(s)}) \theta(f) = 0 .
\end{gather*}
Now set
\[
K(s) = \sup_{t \in [s,s+1)} | f(t) - f(s)|   , \qquad s \geq 0 .
\]
Then
\begin{gather*}
K(s)
  =   \sup_{t \in [s,s+1)} \left| \frac{1}{c(t)} \int_1^{e^t} g^*(u)du -
\frac{1}{c(s)} \int_1^{e^s} g^*(u)du\right| \\
\phantom{K(s)}
  \leq   \frac{1}{c(s)} \int_1^{e^{s+1}} g^*(u)du -
\frac{1}{c(s+1)} \int_1^{e^s} g^*(u)du \\
\phantom{K(s)}
  \leq   \frac{1}{c(s)} \int_{e^s}^{e^{s+1}} g^*(u)du + \left(1-\frac{c(s)}{c(s+1)}\right)
\frac{1}{c(s)} \int_1^{e^s} g^*(u)du \\
\phantom{K(s)}
  \leq   k(s) + \left|1-\frac{c(s)}{c(s+1)}\right|f(s) .
\end{gather*}
Hence
\[
\theta(K) = 0 .
\]
We have
\begin{gather*}
|\theta((1-pr)(f))|   =   |\theta(f(t) - f(\floor{t}))|
  =   |\theta(f(t-1) - f(\floor{t})| \\
\phantom{|\theta((1-pr)(f))|}{}
 \leq   \theta(|f(t-1) - f(\floor{t})|)
  \leq   \theta(K(t-1)) = \theta(K) = 0 .
\end{gather*}
Also{\samepage
\begin{gather*}
|\theta((1-p_cr)(f))|   =   |\theta((t-\floor{t})(f(t) - f(\ceil{t}))
+(1-(t-\floor{t}))(f(t) - f(\floor{t})))| \\
\phantom{|\theta((1-p_cr)(f))|}{}
 \leq   \theta(|f(t) - f(\ceil{t})|) + \theta(|f(t-1) - f(\floor{t})|) \\
\phantom{|\theta((1-p_cr)(f))|}{}
  \leq   \theta(K(t)) + \theta(K(t-1)) = 2 \theta(K) = 0.
\end{gather*}
This demonstrates 1.}

Similarly
\begin{gather*}
|(1-E)(f)(s)|   =   \left|f(s) - \int_{s}^{s+1} f(t)dt\right|
  \leq  \sup_{t \in [s,s+1)}|f(t) - f(s)| = K(s).
\end{gather*}
Hence
\[
|\theta((1-E)(f))| \leq \theta(|(1-E)(f)|) \leq \theta(K) = 0 .
\]
This demonstrates 2.
\end{proof}

\begin{lemma} \label{lemma:app:copyLSS}
For every $\xi \in BL_{(c)}[0,\infty)$ $($resp.~$CBL_{(c)}[0,\infty)$, $B(C)_{(c)}[0,\infty))$ there exists $\xi' \in BL$ $($resp.~$CBL, B(C))$
$($and for every $\xi' \in BL$ $($resp.~$CBL, B(C))$ there exists $\xi \in BL_{(c)}[0,\infty)$ $($resp.~$CBL_{(c)}[0,\infty)$, $B(C)_{(c)}[0,\infty)))$
such that $(\xi-\xi'r)L^{-1}(\alpha_g)(t) =0$ $\fa\, g \in \mathrm{m}_{1,\infty}$.
\end{lemma}
\begin{proof}
Let $\xi \in BL_{(c)}[0,\infty)$. Set $\xi'=\xi p_{(c)}$.  It is easily verif\/ied
$\xi'$ is a state on $\ell^\infty$ which is translation invariant by
Lemma~\ref{lemma:app:prelim}.1. Moreover,
\[
(\xi-\xi'r)L^{-1}(\alpha_g) = \xi(1-p_{(c)}r)L^{-1}(\alpha_g) = 0
\]
due to Lemma~\ref{lemma:app:prelim2}.  Now let $\xi' \in BL$.
Set $\xi=\xi'rE$.  It is easily verif\/ied
$\xi$ is a state on $L^\infty([0,\infty))$ (resp.~$C_b([0,\infty))$) which is translation invariant for $a \in \NN$ by
Lemma~\ref{lemma:app:prelim}.3.  Let $a = j + k$ where $j \in \NN$
and $k \in (0,1)$.  Then
\[
\int_{n+k}^{n+1+k}
f(s)ds = \int_{n+k}^{n+1}f(s)ds + \int_{n+1}^{n+1+k}f(s)ds
\]
and
\begin{gather*}
\xi'rE(T_k f)   =   \xi'\left(\int_{n+k}^{n+1}f(s)ds\right) + \xi'\left(\int_{n+1}^{n+1+k}f(s)ds\right) \\
\phantom{\xi'rE(T_k f)}{}  =   \xi'\left(\int_{n+k}^{n+1}f(s)ds\right) + \xi'\left(\int_{n}^{n+k}f(s)ds\right)
  =   \xi'\left(\int_{n}^{n+1}f(s)ds\right) = \xi'rE(f) .
\end{gather*}
Hence $\xi'rE$ is translation invariant for all $a > 0$.  Moreover,
\[
(\xi-\xi'r)L^{-1}(\alpha_g) = \xi'r(1-E)L^{-1}(\alpha_g) = 0
\]
due to Lemma~\ref{lemma:app:prelim2}.

Now let $\xi = \gamma C \in CBL_{(c)}[0,\infty)$ for any singular state $\gamma$ on $L^\infty([0,\infty))$ (resp.~$C_b([0,\infty))$).  Then
\[
\xi' = \xi p_{(c)} = \gamma C p_{(c)} = \gamma p_{(c)} C
\]
by Lemma~\ref{lemma:app:prelim}.2 where $\gamma p_{(c)}$ is a singular state
on $\ell^\infty$.  Hence $\xi' \in CBL$.
Similarly, if $\xi \in B(C)_{(c)}[0,\infty)$,
\[
\xi'C = \xi p_{(c)}C = \xi Cp_{(c)} = \xi p_{(c)} = \xi'
\]
by Lemma~\ref{lemma:app:prelim}.2 again.  Hence $\xi' \in B(C)$.
Conversely, if $\xi' = \gamma' C \in CBL$ with $\gamma'$ a singular state on $\ell^\infty$ then
\[
\xi = \xi'rE = \gamma' C rE = \gamma' rE C
\]
by Lemma~\ref{lemma:app:prelim}.4 where $\gamma' rE$ is a singular state on $L^\infty([0,\infty))$ (resp.~$C_b([0,\infty))$).
Similarly, if $\xi' \in B(C)$,
\[
\xi C = \xi'rEC = \xi'CrE = \xi'rE = \xi
\]
by Lemma~\ref{lemma:app:prelim}.4 again. Hence $\xi \in B(C)_{(c)}[0,\infty)$.
\end{proof}

\begin{corollary}
We have the equalities as in \eqref{eq:dtr1}, \eqref{eq:dtr2}
and \eqref{eq:dtr3}.
\end{corollary}
\begin{proof}
We note $g = \sum_{n=1}^\infty \mu_n(T) \chi_{[n,n+1)} \in \mathrm{m}_{1,\infty}$
for all $T \in \mathcal{M}_{1,\infty}$.
The result for \eqref{eq:dtr1} follows from the diagram
\begin{equation} \label{fig:app:simplify:1}
\begin{array}{ccccccc}
%\mathcal{V}_{DL}
& %\stackrel{2.}{=}
& \mathcal{V}_{DL[1,\infty)}
& & & & \mathcal{V}_{DL_c[1,\infty)} \\[4pt]
& & L \shortparallel & & & & L \shortparallel \\
& & \mathcal{W}_{BL[0,\infty)} & \stackrel{\ref{lemma:app:copyLSS}}{=} & \mathcal{W}_{BL}
& \stackrel{\ref{lemma:app:copyLSS}}{=} & \mathcal{W}_{BL_c[0,\infty)}
\end{array}				
\end{equation}
where $L$ denotes the isomorphism between
$\mathcal{W}_{BL[1,\infty)}$ (resp.\ $\mathcal{W}_{BL_c[1,\infty)}$)
and $\mathcal{V}_{DL[1,\infty)}$ (resp.\ $\mathcal{V}_{DL_c[1,\infty)}$)
induced by the map $L^{-1}$ \eqref{eq:app:L},
and $\stackrel{\ref{lemma:app:copyLSS}}{=}$ denotes the statement of
Lemma~\ref{lemma:app:copyLSS}.
The result for \eqref{eq:dtr2}
and \eqref{eq:dtr3} follows from an identical diagram.
\end{proof}

\subsection*{Acknowledgements}

The authors thank N.~Kalton for many discussions concerning Dixmier traces and symmetric functionals.

%\pdfbookmark[1]{References}{ref}

\LastPageEnding

\end{document}